\documentclass[draft]{Dekker}
\usepackage{amsmath,amsthm,amscd,amssymb}

\begin{document}

\long\def\symbolfootnote[#1]#2{\begingroup%
\def\thefootnote{\fnsymbol{footnote}}\footnote[#1]{#2}\endgroup}

\newcommand{\Q}{\nabla}
\newcommand{\Div}{\textup{div}}
\newcommand{\Grad}{\textup{grad}}
\newcommand{\Ess}{\textup{ess sup}}
\newcommand{\Min}{\textup{Min}}
\newcommand{\D}{\Delta}
\newcommand{\ka}{\kappa}
\newcommand{\la}{\lambda}
\newcommand{\Om}{\Omega}
\newcommand{\rz}{\mathbb{R}}
\newcommand{\sz}{\mathbb{S}}
\newcommand{\hz}{\mathbb{H}}
\newcommand{\cz}{\mathbb{C}}
\newtheorem{thm}{Theorem}
\newtheorem{cor}[thm]{Corollary}
\newtheorem{lem}[thm]{Lemma}
\newtheorem{prop}[thm]{Proposition}

\setcounter{page}{001}
\jname{
xxx}
\jvol{xx}
\jissue{x}
\jyear{xxxx}
\webslug{xxx}
\cpright{xxx}

\title[Harrell-Stubbe Inequalities]{
On Harrell-Stubbe Type Inequalities for the
Discrete Spectrum of a Self-Adjoint Operator}

\author[Ashbaugh and Hermi]{Mark S.\ Ashbaugh$^{\dagger}$
\symbolfootnote[2]{Partially supported by National Science
Foundation (USA) grant DMS--9870156.} ,
Lotfi Hermi$^{\ddagger}$ \symbolfootnote[3]{Corresponding author.  Part of this
paper appeared in L.H.'s 1999 Ph.D.\ Dissertation.}
\affiliation{$^{\dagger}$ Department of Mathematics\\  University of Missouri\\
Columbia, MO \, 65211-4100, USA\\
Email: mark@math.missouri.edu
\\
$^{\ddagger}$ Department of Mathematics\\  University of Arizona\\
Tucson, AZ \, 85721, USA\\
Email: hermi@math.arizona.edu}}

\abstract{We produce a new proof and extend results by Harrell and
Stubbe for the discrete spectrum of a self-adjoint operator.  An
abstract approach--based on commutator algebra, the Rayleigh-Ritz
principle, and an ``optimal'' usage of the Cauchy-Schwarz
inequality--is used to produce ``parameter-free'',
``projection-free'' versions of their theorems.  We also analyze
the strength of the various inequalities that ensue.  The results
contain classical bounds for the eigenvalues. Extensions of a
variety of inequalities \`a la Harrell-Stubbe are illustrated for
both geometric and physical problems.}

\keywords{eigenvalues of the Laplacian; Dirichlet eigenvalue problem for domains
in Euclidean space; eigenvalues of elliptic operators; Payne--P\'{o}lya-- Weinberger
inequality; Hile--Protter inequality; H.\ C.\ Yang inequality; Harrell-Stubbe inequalities;
universal eigenvalue estimates; reverse Chebyshev inequality.}

\maketitle
%\today

\section{Introduction}\label{intro}

In this paper, we continue our work started in \cite{AH2}.  A semibounded operator
modeled after the Dirichlet Laplacian on a bounded domain $\Om \subset \rz^n$
(or a Schr\"odinger operator with magnetic potential) is given.  We provide universal
bounds for its eigenvalues.  These are estimates for the eigenvalues that do not
involve domain dependencies \cite{Pr3} (see also \cite{A3}, \cite{A2}).  This is a
problem related to a classical result of Payne, P\'olya, and Weinberger \cite{PPW1},
\cite{PPW2} (abbreviated as PPW) for the eigenvalues $0< \la_1 < \la_2 \le \la_3
\le \cdots$ (multiplicities included) of the fixed membrane problem

\begin{align}
 - \Delta u & = \lambda u
\quad \text{in } \quad \Omega, \nonumber \\
u & = 0 \quad \, \text{ on } \quad \partial \Omega.
\label{eq:intro.2}
\end{align}

We provide, based on the Rayleigh-Ritz principle and trial functions, alternative
proofs and extensions of recent results which were obtained by Harrell and Stubbe
\cite{HS}.  Our main divergence from their method is the use of the ``optimal''
Cauchy-Schwarz inequality exploited in \cite{AH2} (see also \cite{A3}, \cite{A2},
\cite{AB1}, and \cite{Y}) and the fact that we employ the Rayleigh-Ritz inequality and
not algebraic identities (see also \cite{LP} for yet another alternative).

We also consider the inequalities
\begin{equation}
\sum_{i=1}^m (\la_{m+1} - \la_i)^p \leq
\dfrac{2 p}{n} \sum_{i=1}^m  \la_i (\la_{m+1} - \la_i)^{p-1}
\quad \text{ for } \quad  p \ge 2
\label{eq:intro.5}
\end{equation}
(see ineq.\ (14) in Theorem9, p.\ 1805 of \cite{HS}) and
\begin{equation}
\sum_{i=1}^m (\la_{m+1} - \la_i)^p \leq \dfrac{4}{n} \sum_{i=1}^m
\la_i (\la_{m+1} - \la_i)^{p-1} \quad \text{ for } \quad   p \le 2
\label{eq:intro.6}
\end{equation}
(see ineq.\ (11) in Theorem 5, p.\ 1801 of \cite{HS}), which stem from two different
considerations in Harrell and Stubbe's work.

The classical Hile-Protter \cite{HP} and H.\ C.\ Yang \cite{Y} inequalities appear as
special cases of (\ref{eq:intro.6}) for $p=0$ and $2$ respectively.  In this paper, we
will in fact show that (\ref{eq:intro.6}) improves monotonically for $0 \le p \le 2$
(Theorem \ref{thm:4.1}).  The classical PPW inequality

\begin{equation}
\la_{m+1} - \la_{m} \leq  \dfrac{4}{n} \, \frac{\sum_{i=1}^{m} \la_i}{m},
\label{eq:intro.1}
\end{equation}
is obviously weaker than the $p=1$ case of (\ref{eq:intro.6}).  In fact, it is weaker
than the $p=0$ case of (\ref{eq:intro.6}) (see \cite{HP}), which is also easy to see.

In the literature, the case $p=1$ is referred to as the {\em ``Yang 2''} bound (see
\cite{A3}, \cite{A2}).  It is, of course, explicitly given by
\begin{equation}
\la_{m+1}  \leq  \left(1+ \dfrac{4}{n}\right) \, \dfrac{\sum_{i=1}^{m} \la_i}{m},
\label{eq:intro.yang2}
\end{equation}
The case $p=0$ (herein referred to as HP) reads, explicitly,
\begin{equation}
\dfrac{m n}{4} \leq \sum_{i=1}^m  \dfrac{\la_i}{\la_{m+1} - \la_i}.
\label{eq:intro.hp}
\end{equation}

The general framework for this paper provides extensions \`a la
Harrell-Stubbe for various geometric and physical problems.  In fact
Harrell-Stubbe type inequalities are valid for all the situations for
which H.~C.~Yang-style improvements have been proved in \cite{AH2} and
illustrated in \cite{AH4}.

In this paper, an analysis of (\ref{eq:intro.5}) is also provided.
It is proved that the case $p=2$ (i.e., the H.~C.~Yang inequality,
also referred to as {\em ``Yang 1''}; see \cite{A2}, \cite{A3}) is
the strongest for $p\ge 2$ (Theorem \ref{thm:4.2}).

\section{General Framework}  \label{sec:setup}

Our setting is that of \cite{AH2}.  We provide an ``algebraized''
version of the membrane problem described in the introduction.
Such a scheme follows a line of thought first adopted by Harrell
and Davies (see \cite{Ha}, \cite{Michel}) in 1988.  This
abstraction has the advantage of unifying many results for gaps of
eigenvalues of subdomains of Riemannian manifolds and a variety of
geometric and physical situations.  This point of view was
advocated by Harrell and Michel \cite{HM2}, \cite{HM1},
\cite{Michel}, Harrell and Stubbe \cite{HS}, Hook \cite{Ho2},
Levitin and Parnovski \cite{LP}, and Ashbaugh and Hermi
\cite{AH2}.  This point of view provides improvements \`a la
H.~C.~Yang of results in \cite{Cheng}, \cite{HM2}, \cite{HM1},
\cite{Ho2}, \cite{Lee}, \cite{Leung}, \cite{Li}, and \cite{YY} as
described in \cite{AH2} and \cite{AH4}. See also \cite{HarCPDE}
and \cite{HarSoufiIlias} where further applications and
generalizations to new settings are presented.

A complex Hilbert space $\mathcal{H}$ with inner product $\langle \ , \ \rangle$ is
given.  $\langle \ , \ \rangle$ is taken to be linear in its first argument, conjugate linear
in its second.  We let $A: \mathcal{D} \subset \mathcal{H} \rightarrow \mathcal{H} $ be
a self-adjoint operator defined on a dense domain $\mathcal{D}$ which is
semibounded below and has a discrete spectrum $\lambda_1 \leq \lambda_2
\leq \lambda_3 \leq \cdots.$  Let $\{B_k : A(\mathcal{D}) \rightarrow
\mathcal{H} \}_{k=1}^{N}$ be a collection of symmetric operators which leave
$\mathcal{D}$ invariant, and let $\{u_i\}_{i=1}^{\infty}$ be the normalized
eigenvectors of $A$, $u_i$ corresponding to $\la_i$.  This family of eigenvectors
is further assumed to be an orthonormal basis for $\mathcal{H}$.  The commutator
of two operators, $[A,B]$, is defined by $[A,B] = AB- BA$, and
$\|u \| = \sqrt{\langle u,u \rangle}$.

As in \cite{AH2}, we define
\begin{align}
\rho_i = \sum_{k=1}^{N} \langle [A,B_k] u_i, B_k u_i\rangle
\label{eq:1.1}
\end{align}
and
\begin{align}
\Lambda_i = \sum_{k=1}^{N} \| [A,B_k] u_i \|^2.
\label{eq:1.2}
\end{align}

In \cite{AH2}, we have shown that the classical inequalities of PPW, HP, and
H.\ C.\ Yang follow from the same general set-up and the following theorem.

\begin{thm}
 \label{thm:1.1}
%Theorem 1
The eigenvalues $\la_i$ of the operator $A$ satisfy the inequalities
\begin{align}
\sum_{i=1}^m \rho_i
\leq \dfrac{ \sum_{i=1}^m \Lambda_i } {\la_{m+1} - \la_m},
\label{eq:1.3}
\end{align}
\begin{align}
\sum_{i=1}^m \rho_i
\leq \sum_{i=1}^m \dfrac{\Lambda_i } {\la_{m+1} - \la_i},
\label{eq:1.4}
\end{align}
and
\begin{align}
\sum_{i=1}^m (\la_{m+1} - \la_i)^2 \rho_i
\leq \sum_{i=1}^m (\la_{m+1} - \la_i) \Lambda_i.
\label{eq:1.5}
\end{align}
\end{thm}

These give abstract versions of the PPW, HP, and Yang inequalities,
respectively, and even at this level, (\ref{eq:1.5}) is stronger than
(\ref{eq:1.3}) and (\ref{eq:1.4}), and (\ref{eq:1.4}) is stronger than
(\ref{eq:1.3}).

\section{Extending The Work of Harrell and Stubbe}  \label{sec:3}

Based solely on the ``traditional'' tools (the Rayleigh-Ritz principle, simple
trial functions, the Cauchy-Schwarz inequality, $\dots$), in this section
we provide alternative proofs and generalizations of the results of Harrell and
Stubbe \cite{HS}.  In their work, they wanted to understand the nature of Yang's
inequalities \cite{Y}.  Our proofs tie in with the abstract commutator approach
used by various authors \cite{HM2}, \cite{HM1}, \cite{HS}, in their work on
geometric bounds for eigenvalues of elliptic operators.  The proofs provide further
insight into the extensions in \cite{HS} (explaining, for example, what terms
are being dropped in arriving at their inequalities).  A separate section (Section
\ref{sec:bounds}) is dedicated to comparing the bounds obtained from the
approach given in this section to those of the works of Hile-Protter and H.\ C.\ Yang.
Another section (Section \ref{applications}) provides illustrations of various
extensions of known bounds for geometric and physical problems; for more in
this direction see \cite{AH4}.

\begin{thm}
\label{thm:3.1}
%Theorem 2
Let the function $g(\la)$ be nonnegative and nondecreasing on the eigenvalues
$\{\la_i\}_{i=1}^m$ of $A.$  Then the eigenvalues $\{\la_i\}_{i=1}^{m+1}$
of $A$ satisfy the inequality
\begin{align}
\sum_{i=1}^{m} (\la_{m+1} - \la_i)^2 g(\la_i) \rho_i \leq
\sum_{i=1}^{m} (\la_{m+1} - \la_i) g(\la_i) \Lambda_i.
\label{eq:3.1}
\end{align}
\end{thm}

\smallskip \noindent {\em Note.}  It is enough that $g$ be a nonnegative
and nondecreasing function defined on $(0,\la_{m+1}),$ as will typically
be the case in applications.

\smallskip \noindent {\em Proof.}
As in \cite{AH2}, we start with the Rayleigh-Ritz inequality
\begin{align}
\la_{m+1} \leq \dfrac{\langle A \phi, \phi \rangle}
{\langle \phi, \phi \rangle}
\label{eq:2.1}
\end{align}
and the test function
\begin{align}
\phi_i = B u_i - \sum_{j=1}^m a_{ij} u_j,
\label{eq:2.3}
\end{align}
where $B$ is one of the $B_k$'s, $k=1, ..., N$.  The orthogonality condition
\[{\langle \phi, u_j \rangle} = 0\] for $j=1, 2, \cdots, m$ makes $a_{ij} =
\langle B u_i, u_j \rangle.$  The symmetry of $B$ makes $a_{ji} =
\overline{a_{ij}}.$  As in \cite{AH2}, (\ref{eq:2.1}) reduces to
\begin{align}
\la_{m+1} - \la_i \leq \dfrac{\langle [A,B] u_i, \phi_i \rangle}
{\langle \phi_i, \phi_i \rangle}.
\label{eq:2.6}
\end{align}
The calculations in \cite{AH2} yield
\begin{align}
\langle [A,B] u_i, \phi_i \rangle =
\langle [A,B] u_i, B u_i \rangle - \sum_{j=1}^m (\la_j - \la_i) \ |a_{ij}|^2.
\label{eq:2.8}
\end{align}
By (\ref{eq:2.6}), $\langle [A,B] u_i, \phi_i \rangle \ge 0$.  Thus, by the
``optimal'' Cauchy-Schwarz inequality (see Lemma 3.1 of \cite{AH2}),
\begin{align}
\dfrac{\langle [A,B] u_i, \phi_i \rangle}{\langle \phi_i, \phi_i \rangle}
\leq \dfrac{\| [A,B] u_i \|^2 -\sum_{j=1}^m (\la_j- \la_i)^2 \ |a_{ij}|^2
}{\langle [A,B] u_i, \phi_i \rangle}.
\label{eq:2.11}
\end{align}
We then obtain
\begin{align}
(\la_{m+1} - \la_i) \Big( \langle [A,B] u_i, B u_i \rangle
- &\sum_{j=1}^m (\la_j - \la_i) |a_{ij}|^2 \Big) \notag \\
&\leq \| [A,B] u_i \|^2 -\sum_{j=1}^m (\la_j - \la_i)^2 |a_{ij}|^2,
\label{eq:2.12}
\end{align}
or, upon combining the sums involving $|a_{ij}|^2,$
\begin{align}
(\la_{m+1} - \la_i) \langle [A,B] u_i, B u_i \rangle \notag \\
&\leq \| [A,B] u_i \|^2 -\sum_{j=1}^m (\la_i - \la_j)
(\la_{m+1} - \la_j)|a_{ij}|^2,
\label{eq:2.12a}
\end{align}
Since $B$ is one of the $B_k$'s, $a_{ij} \equiv a_{ij}^k.$  Let
\begin{align}
A_{ij} \equiv \sum_{k=1}^N |a_{ij}^k|^2.
\label{eq:2.120}
\end{align}
Hence $A_{ji}=A_{ij} \ge 0.$  Replacing $B$ by $B_k$ in (\ref{eq:2.12a}), summing
over $k$ for $1 \le k \le N$, and incorporating the definitions of $\rho_i$, $\Lambda_i$,
and $A_{ij},$ we obtain
%\begin{align}
%(\la_{m+1} - \la_i) \Big( \rho_i - \sum_{j=1}^m (\la_j - \la_i) A_{ij} \Big)
%\leq \Lambda_i -\sum_{j=1}^m (\la_j - \la_i)^2 A_{ij}.
%\label{eq:2.13}
%\end{align}
%We rewrite (\ref{eq:2.13}) as
\begin{align}
(\la_{m+1} - \la_i) \rho_i
\leq \Lambda_i - \sum_{j=1}^m (\la_{m+1} - \la_j) (\la_i - \la_j) A_{ij}.
\label{eq:2.18}
\end{align}

Multiplying both sides by \[(\la_{m+1} - \la_i) \ g(\la_i)\] ($g \ge 0$ is needed here
to preserve the sense of our inequality, and of course $i<m+1$ is assumed) and
summing over $i$, $1 \leq i \leq m$, gives
\begin{align}
\sum_{i=1}^{m} (\la_{m+1} - \la_i)^2 &g(\la_i) \rho_i \leq
\sum_{i=1}^{m} (\la_{m+1} - \la_i) g(\la_i) \Lambda_i \notag \\
&- \sum_{i,j=1}^{m} (\la_{m+1} - \la_i)
(\la_{m+1} - \la_j) A_{ij} (\la_i-\la_j) g(\la_i).
\label{eq:3.2}
\end{align}

If $g$ were constant the double sum in $i$ and $j$ here would vanish due to
antisymmetry (recall that $A_{ij}$ is symmetric), allowing us to conclude that the
theorem holds in this case.  Indeed, it is this case which motivated our choice of
multiplier for (\ref{eq:2.18}).  This is how Yang's main inequality (the $p=2$ case
of (\ref{eq:intro.6})) was proved in \cite{A3} and \cite{A2}.  For more general $g,$
we can use the notions of {\it similarly} (resp., {\it oppositely}) {\it ordered} (see
\cite{HLP}, pp.\ 43, 261-262) to impose a sign on the double sum; in particular, it
transpires that the double sum is nonnegative if $\{\la_i\}_{i=1}^m$ and
$\{g(\la_i)\}_{i=1}^m$ are similarly ordered, or, what amounts to the same thing
here, if $\{g(\la_i)\}_{i=1}^m$ is a nondecreasing sequence.  To see this, we
rewrite the double sum with $i$ and $j$ interchanged and average the two
expressions, giving
\begin{align}
\sum_{i=1}^{m}(\la_{m+1} - \la_i)^2 g(\la_i) \rho_i &\leq
\sum_{i=1}^{m} (\la_{m+1} - \la_i) g(\la_i) \Lambda_i \notag \\
&- \frac{1}{2}
\sum_{i,j=1}^{m} (\la_{m+1} - \la_i) (\la_{m+1} - \la_j) A_{ij} \notag \\
&\times  (\la_i-\la_j) ( g(\la_i)-g(\la_j)).
\label{eq:3.3}
\end{align}
The factor $(\la_i-\la_j)(g(\la_i)-g(\la_j))$ and the nonnegativity of the rest, shows
that the double sum will be nonnegative whenever $\{g(\la_i)\}_{i=1}^m$ is
nondecreasing, and, since the double sum is preceded by a minus sign, its
contribution to the right-hand side of the inequality is nonpositive, yielding the
desired conclusion, i.e., inequality (\ref{eq:3.1}).
%Thus, under the hypotheses of the theorem, the second term on the
%right-hand side of (\ref{eq:3.3}) is nonpositive and the conclusion of
%the theorem follows.
\hfill $\Box$

\vskip .5 cm

\noindent {\em Remarks.} 1.  If one assumes that $g$ is nondecreasing
and $C^1$ (or just differentiable) on the positive half-axis, then by the
mean value theorem
\[g(\la_i)-g(\la_j) = g'(\xi_{ij}) (\la_i-\la_j)\]
for some $\xi_{ij} > 0$ where $g'(\xi_{ij}) \ge 0.$  Therefore the double
sum giving the second term on the right-hand side of (\ref{eq:3.3}) is
nonnegative, and since it is subtracted, the statement of the theorem
follows.
%\hfill $\Box$
%\noindent {\em Remarks.} 1.  The condition that $h(\la)$ be a
%$C^1$-function can be relaxed.  Indeed, since $h(\la)$ is a nonincreasing
%function, the sequences $\la_i$ and $b_i\equiv h(\la_{m+1} - \la_i)$ are
%{\em similarly ordered} (see p.~43 of \cite{HLP}), viz.\,
%\[ (\la_i-\la_j) \  \left( h(\la_{m+1} - \la_i)-h(\la_{m+1} - \la_j) \right)
%\ge 0.\]
%This implies that the dropped term in (\ref{eq:3.3}) is nonpositive.
%To simplify the notation, let $\omega_i=\la_{m+1} - \la_i,
%b_i=h(\la_{m+1} - \la_i)$.
%%This defines b_i, and also corrects an index in the definition of \omega_i.
%The dropped term then takes the form
%\begin{align}
%\frac{1}{2}
%\sum_{i,j=1}^{m} \omega_i \ \omega_j \ A_{ij}
%(\la_j-\la_i) \  (b_i-b_j) &= - \frac{1}{2} \sum_{i,j=1}^{m} \omega_i \ \omega_j \ A_{ij}
%(\la_i-\la_j) \  (b_i-b_j) \notag \\
%&\le - \frac{1}{2} \ \underset{1\le i,j \le m}{\Min} A_{ij} \ \sum_{i,j=1}^{m} \omega_i \omega_j
%(\la_i-\la_j) \ (b_i-b_j) \notag \\
%&\le 0 \notag
%\end{align}

\noindent 2.  If we make the replacement $g(\la)=h(\la_{m+1}-\la),$ then the
hypotheses on $h$ would be that $h$ is nonnegative and nonincreasing on the
sequence $\{\la_{m+1}-\la_i\}_{i=1}^m,$ or, perhaps a little more naturally,
that $h$ is nonnegative and nonincreasing on $(0, \la_{m+1}).$  The inequality
in Theorem \ref{thm:3.1},
%Theorem 2
when written in terms of $h,$ becomes
\begin{align}
\sum_{i=1}^{m} (\la_{m+1}-\la_i)^2 h(\la_{m+1} - \la_i) \rho_i \leq
\sum_{i=1}^{m} (\la_{m+1}-\la_i) h(\la_{m+1} - \la_i) \Lambda_i.
\label{eq:3.4}
\end{align}

\noindent 3.  Setting $f(\la)= (\la_{m+1} - \la)^2 g(\la),$ or,
equivalently, in the notation of Remark 2,
$f(\la)= (\la_{m+1} - \la)^2 h(\la_{m+1} - \la),$ for
$0 < \la \leq \la_{m+1},$
(\ref{eq:3.1}) can be written as
\begin{align}
\sum_{i=1}^{m} f(\la_i) \rho_i \leq
\sum_{i=1}^{m} \dfrac{f(\la_i)}{\la_{m+1} - \la_i} \Lambda_i.
\label{eq:3.4a}
\end{align}
This is the statement of Theorem 5 in \cite{HS} (when one specializes to their
setting, which leads to $\rho_i=N, \Lambda_i=4\la_i,$ as in our Corollary
\ref{cor:3.3}
%Corollary 4
below).  The condition that the function $f(\la)(\la_{m+1}- \la)^{-2}$ (in their case)
be nondecreasing is equivalent to the statement of Theorem \ref{thm:3.1} which
seems more natural to the problem.  The H.\ C.\ Yang type inequality (\ref{eq:1.5})
obtains when $f(\la)= (\la_{m+1} - \la)^2$, i.e., when $g(\la) \equiv 1$ (or,
equivalently, when $h(\la) \equiv 1$).  As noted earlier, the second term on the
right-hand side of (\ref{eq:3.3}) is identically zero in this case.

\begin{cor}
\label{cor:3.2}
%Corollary 3
Let $p\le 2$.  Then
%Let $0 \le p \le 2$.  Then
\begin{align}
\sum_{i=1}^m (\la_{m+1} - \la_i)^p \  \rho_i
\leq \sum_{i=1}^m (\la_{m+1} - \la_i)^{p-1} \ \Lambda_i.
\label{eq:3.142}
\end{align}
\end{cor}
\smallskip \noindent {\em Proof.}
We make the choice $g(\la)=(\la_{m+1}-\la)^{p-2}$ (or equivalently $h(\la)= \la^{p-2}$
if applying Remark 2) for $\la \ge 0$ in Theorem \ref{thm:3.1}.
\hfill $\Box$

\vskip 0.5 cm

\noindent {\em Remark.} We note that (\ref{eq:1.5}) and (\ref{eq:1.4}) are particular
instances of this corollary (for $p=2$ and $p=0$, respectively) while (\ref{eq:1.3}) is
a weaker result obtained from (\ref{eq:1.4}) by replacing $\la_{m+1}-\la_i$ by
$\la_{m+1}-\la_m$.

\begin{cor}
\label{cor:3.3}
%Corollary 4
Suppose $A= - \sum_{k=1}^{N} T_k^2$ where the $T_k$'s are skew- symmetric with
domains $\mathcal{D}(T_k)$ such that $\mathcal{D} = \mathcal{D}(A) \subset
\mathcal{D} (T_k)$ and $T_k (\mathcal{D}) \subset \mathcal{D} (T_k)$ and suppose
that $[T_\ell, B_k] u = \delta_{\ell k} u.$  Then $\rho_i = N$, $\Lambda_i = 4 \la_i$, and
for $p\le 2$
%for $0 \le p \le2$
\begin{align}
\sum_{i=1}^m (\la_{m+1} - \la_i)^p
\leq \dfrac{4}{N} \sum_{i=1}^m (\la_{m+1} - \la_i)^{p-1} \ \la_i.
\label{eq:3.6}
\end{align}
\end{cor}
\smallskip \noindent {\em Proof.}  The details of the calculations of $\rho_i$ and
$\Lambda_i$ in this case are provided by Lemma 2.3 and Corollary 2.4 of
\cite{AH2}.
\hfill $\Box$

\vskip 0.5 cm

\noindent {\em Remark.}  Inequality (\ref{eq:intro.6}) is a particular case of this
corollary with $A=-\D$, $T_k=\frac{\partial}{\partial x_k}$, $B_k=x_k$, and $N=n$,
the spatial dimension.

For a symmetric operator $C$ and $\alpha \in \rz$, $C \geq \alpha$ if
$\langle C u, u \rangle \geq \alpha \langle u, u \rangle$ for all vectors
$u\in \mathcal{D} (C)$.  Moreover, $A \geq B$ for symmetric operators $A$ and $B$
if $\mathcal{D} (B) \subset \mathcal{D}(A)$ and $A-B \ge 0$ on $\mathcal{D} (B)$.

\begin{cor}
\label{cor:3.4}
%Corollary 5
Suppose there exist $\gamma, \beta$ such that
\begin{align}
0 < \gamma \leq  [B_k, [A, B_k]]
\label{eq:3.70}
\end{align}
and
\begin{align}
- \sum_{k=1}^N [A, B_k]^2 \leq \beta A.
\label{eq:3.80}
\end{align}
Then, for $p \le 2$,
%Then, for $0 \le p \le 2$,
\begin{align}
\sum_{i=1}^m (\la_{m+1} - \la_i)^p
\leq \dfrac{2 \beta}{\gamma N} \sum_{i=1}^m (\la_{m+1} - \la_i)^{p-1} \ \la_i.
\label{eq:3.91}
\end{align}
\end{cor}
\smallskip \noindent {\em Proof.}
As observed in \cite{AH2} (see Theorem~2.5) the conditions of this corollary yield
immediately $\rho_i \ge \dfrac{1}{2} \gamma N$ and $\Lambda_i \le \beta \la_i.$
Substituting these inequalities into Corollary \ref{cor:3.2}, we obtain the desired
result.
\hfill $\Box$

\smallskip

We now deal with a second set of inequalities treated by Harrell-Stubbe in \cite{HS}.
We adopt their definition:  A real function $f(x)$ is said to satisfy condition (H1) if there
exists a function $r(x)$ such that
\begin{align}
\text{(H1)} \qquad \qquad  \dfrac{f(x)- f(y)}{x-y} \geq \dfrac{r(x)+ r(y)}{2}. \notag
\end{align}
As an example, a function whose derivative $f'$ is concave satisfies this condition.

\begin{lem} \label{lemma:3.2}
Suppose a $C^1$ function $f$ is such that its derivative $f'$ is concave.  Then $f$
satisfies (H1) with $r(x)=f'(x)$.
\end{lem}
\smallskip \noindent {\em Proof.}
For each $\xi$ between $x$ and $y$, $\exists ! \ \mu \in [0,1]$ such that
$\xi= \mu \ x + (1-\mu) \ y$.
Without loss of generality, we can assume $x \geq y$.  Since $f'$ is concave
\begin{align}
f'(\xi) \geq \mu \ f'(x) + (1-\mu) \ f'(y). \notag
\end{align}
Integrating over $\xi$ from $y$ to $x$ yields (on the right we integrate in $\mu$ from
$0$ to $1$ noting that $d \xi = (x-y) d \mu$)
\[f(x) - f(y) \geq \frac{1}{2} (x-y) ( f'(x) + f'(y)),
\]
and the lemma is immediate.
\hfill $\Box$

\begin{thm} \label{thm:3.3}
Let $f(x)$ be an (H1) function for some $r(x)$.  Then, we have
\begin{equation}
\sum_{i=1}^{m} f(\la_i) \rho_i \leq
- \frac{1}{2} \sum_{i=1}^{m} r(\la_i) \Lambda_i  + \mathcal{R},
\label{eq:3.19}
\end{equation}
where
\begin{equation}
\mathcal{R}= \sum_{k=1}^{N} \sum_{i=1}^m \sum_{j=m+1}^{\infty}
|\langle [A, B_k] u_i, u_j \rangle|^2
\left( \dfrac{f(\la_i)}{\la_{m+1} - \la_i} + \dfrac{1}{2} r(\la_i) \right).
\label{eq:3.51}
\end{equation}
\end{thm}

\noindent {\em Remark.}  We recall that, since $\{u_i\}_{i=1}^{\infty}$ is a basis for
$\mathcal{H}$, one may write
\begin{align}
[A, B_k] u_i = \sum_{j=1}^{\infty} \langle [A, B_k] u_i, u_j \rangle u_j.
\label{eq:3.60}
\end{align}
Furthermore, we have
\begin{align}
\| [A,B_k] u_i \|^2= \sum_{j=1}^{\infty} |\langle [A, B_k] u_i, u_j \rangle |^2.
\label{eq:3.7}
\end{align}
This makes
\begin{align}
\sum_{j=m+1}^{\infty}
|\langle [A, B_k] u_i, u_j \rangle|^2 <\infty
\label{eq:3.71}
\end{align}
for each $i=1, \cdots, m.$  Thus the expression for $\mathcal{R}$ given above is
well-defined.

\smallskip \noindent {\em Proof.}
With the substitution $f(\la) = (\la_{m+1} - \la)^2 g(\la)$ and, {\em a priori}, no
conditions on the function $g(\la)$, calculations down to (\ref{eq:3.3}) can be carried
out as above.

Recalling the definitions of $\Lambda_i$ in (\ref{eq:1.2}) and that of $A_{ij}$ in
(\ref{eq:2.120}), we rewrite (\ref{eq:3.3}) in the form
\begin{align}
\sum_{i=1}^{m} f(\la_i) &\rho_i \leq
\sum_{k=1}^{N} \sum_{i=1}^{m} \dfrac{f(\la_i)}{\la_{m+1} - \la_i}
\| [A,B_k] u_i \|^2 \notag \\
&- \frac{1}{2}
\sum_{k=1}^{N} \sum_{i,j=1}^{m} (\la_{m+1} - \la_i) (\la_{m+1} - \la_j)  (\la_i-\la_j)
|\langle B_k u_i, u_j \rangle |^2  \notag \\
&\times \Big( \dfrac{f(\la_i)}{(\la_{m+1} - \la_i)^2}- \dfrac{f(\la_j)}
{(\la_{m+1} - \la_j)^2} \Big).
\label{eq:3.8}
\end{align}
The gap formula
$\langle [A, B_k] u_i, u_j \rangle = (\la_j - \la_i) \langle B_k u_i, u_j \rangle$ gives
\begin{align}
\sum_{i=1}^{m} f(\la_i) &\rho_i \leq
\sum_{k=1}^{N} \sum_{i=1}^{m} \dfrac{f(\la_i)}{\la_{m+1} - \la_i}
\| [A,B_k] u_i \|^2 \notag \\
&- \frac{1}{2}
\sum_{k=1}^{N} \sum_{i,j=1}^{m} (\la_{m+1} - \la_i) (\la_{m+1} - \la_j)
\dfrac{|\langle [A, B_k] u_i, u_j \rangle |^2}{\la_i-\la_j}  \notag \\
&\times \Big(\dfrac{f(\la_i)}{(\la_{m+1} - \la_i)^2}- \dfrac{f(\la_j)}{(\la_{m+1} - \la_j)^2} \Big).
\label{eq:3.9}
\end{align}
(If $\la_i$ ever equals $\la_j$ here, one should interpret the term(s) in which this occurs as
$0$ by using the gap formula in reverse.)  As noted in \cite{HS}, the second term on the
right-hand side can be reduced to
\begin{align}
-\frac{1}{2}
\sum_{k=1}^{N} \sum_{i,j=1}^{m} |\langle [A, B_k] u_i, u_j \rangle |^2
&\Big( \dfrac{f(\la_i) - f(\la_j)}{\la_i - \la_j} + \dfrac{f(\la_j)}{\la_{m+1} - \la_j} +
\dfrac{f(\la_i)}{\la_{m+1} - \la_i} \Big).
\label{eq:3.10}
\end{align}
Symmetry of this expression in $i$ and $j$ reduces (\ref{eq:3.9}) to
\begin{align}
\sum_{i=1}^{m} f(\la_i) &\rho_i \leq
\sum_{k=1}^{N} \sum_{i=1}^{m} \dfrac{f(\la_i)}{\la_{m+1} - \la_i} \| [A,B_k] u_i \|^2 \notag \\
&- \frac{1}{2}
\sum_{k=1}^{N} \sum_{i,j=1}^{m} |\langle [A, B_k] u_i, u_j \rangle |^2
\dfrac{f(\la_i) - f(\la_j)}{\la_i - \la_j} \notag \\
&- \sum_{k=1}^{N} \sum_{i,j=1}^{m} |\langle [A, B_k] u_i, u_j \rangle |^2
\dfrac{f(\la_i)}{\la_{m+1} - \la_i},
\label{eq:3.11}
\end{align}
i.e.,
\begin{align}
\sum_{i=1}^{m} f(\la_i) &\rho_i \leq
\sum_{k=1}^{N} \sum_{i=1}^{m} \dfrac{f(\la_i)}{\la_{m+1} - \la_i}
\Big( \| [A,B_k] u_i \|^2 - \sum_{j=1}^m |\langle [A, B_k] u_i, u_j \rangle |^2  \Big) \notag \\
&- \frac{1}{2}
\sum_{k=1}^{N} \sum_{i,j=1}^{m} |\langle [A, B_k] u_i, u_j \rangle |^2
\dfrac{f(\la_i) - f(\la_j)}{\la_i - \la_j}.
\label{eq:3.12}
\end{align}
Since $f$ satisfies condition (H1), this reduces to
\begin{align}
\sum_{i=1}^{m} f(\la_i) &\rho_i \leq
\sum_{k=1}^{N} \sum_{i=1}^{m} \dfrac{f(\la_i)}{\la_{m+1} - \la_i}
\Big( \| [A,B_k] u_i \|^2 - \sum_{j=1}^m |\langle [A, B_k] u_i, u_j \rangle |^2  \Big) \notag \\
&- \frac{1}{4}
\sum_{k=1}^{N} \sum_{i,j=1}^{m} |\langle [A, B_k] u_i, u_j \rangle |^2
\big( r(\la_i) + r(\la_j) \big).
\label{eq:3.13}
\end{align}
Symmetry in $i$ and $j$ reduces the second term of the right-hand side to
\begin{align}
- \frac{1}{2}
\sum_{k=1}^{N} \sum_{i,j=1}^{m} r(\la_i) |\langle [A, B_k] u_i, u_j \rangle |^2.
\label{eq:3.14}
\end{align}
This, with identity (\ref{eq:3.7}), gives
\begin{align}
\sum_{i=1}^{m} f(\la_i) \rho_i &\leq
\sum_{k=1}^{N} \sum_{i=1}^{m} \sum_{j=m+1}^{\infty} \dfrac{f(\la_i)}{\la_{m+1} - \la_i}
|\langle [A, B_k] u_i, u_j \rangle |^2 \notag \\
&- \frac{1}{2} \sum_{k=1}^{N} \sum_{i,j=1}^{m} r(\la_i) |\langle [A, B_k] u_i, u_j \rangle |^2.
\label{eq:3.15}
\end{align}
Noting that (\ref{eq:3.7}) gives
\begin{align}
\sum_{j=1}^{m} |\langle [A, B_k] u_i, u_j \rangle |^2 = \| [A,B_k] u_i \|^2 -
\sum_{j=m+1}^{\infty} |\langle [A, B_k] u_i, u_j \rangle |^2,
\label{eq:3.16}
\end{align}
we obtain
\begin{align}
\sum_{i=1}^{m} f(\la_i) &\rho_i \leq
\sum_{k=1}^{N} \sum_{i=1}^{m} \sum_{j=m+1}^{\infty}
\dfrac{f(\la_i)}{\la_{m+1} - \la_i}
|\langle [A, B_k] u_i, u_j \rangle |^2 \notag \\
&+ \dfrac{1}{2} \sum_{k=1}^{N} \sum_{i=1}^{m} \sum_{j=m+1}^{\infty} r(\la_i)
|\langle [A, B_k] u_i, u_j \rangle |^2 \notag \\
&- \frac{1}{2} \sum_{k=1}^{N} \sum_{i=1}^{m} r(\la_i) \| [A,B_k] u_i \|^2,
\label{eq:3.17}
\end{align}
which, upon incorporating the definitions of $\Lambda_i$ and $\mathcal{R}$,
is the statement of the theorem.
\hfill $\Box$

\begin{cor}
\label{cor:3.5}
%Corollary 8
Let $p \geq 2$, then
\begin{align}
\sum_{i=1}^m  (\la_{m+1} - \la_i)^p \rho_i
\leq \dfrac{p}{2} \sum_{i=1}^m (\la_{m+1} - \la_i)^{p-1} \Lambda_i.
\label{eq:3.18}
\end{align}
\end{cor}
\smallskip \noindent {\em Proof.}
For $p\ge 2$, $\la \leq \la_{m+1}$, $f(\la) = (\la_{m+1}- \la)^p$, is
such that
\[f'(\la) = - p (\la_{m+1}- \la)^{p-1}\] is concave, and
\[\dfrac{f(\la)}{\la_{m+1} - \la} + \frac{1}{2} r(\la) = \left(1- \frac{p}{2}\right)
(\la_{m+1}- \la)^{p-1} \leq 0.
\]
Thus $\mathcal{R}\le 0$ and inequality (\ref{eq:3.19}) completes the proof.
\hfill $\Box$

Using the same function $f(\la)$ as in the proof of Corollary \ref{cor:3.5}
and the skew-symmetric operators of Corollary \ref{cor:3.3} yields an
analog of Corollary \ref{cor:3.3} for the case $p \ge 2$:
\begin{cor}
 \label{cor:3.6}
%Corollary 9
Suppose $A= - \sum_{k=1}^{N} T_k^2$ where the $T_k$'s are
skew- symmetric with the same conditions as those of Corollary
\ref{cor:3.3}.  Then
for $p\ge 2$
\begin{align}
\sum_{i=1}^m (\la_{m+1} - \la_i)^p
\leq \dfrac{2p}{N} \sum_{i=1}^m (\la_{m+1} - \la_i)^{p-1} \ \la_i.
\label{eq:3.61}
\end{align}
\end{cor}

Similarly, one obtains a $p \ge 2$ analog of Corollary \ref{cor:3.4}:
%Corollary 5
\begin{cor}
\label{cor:3.7}
%Corollary 10
Let $p \geq 2$, and suppose there exist $\gamma$, $\beta$ such that
the conditions of Corollary \ref{cor:3.4} are satisfied, then
%Corollary 5
\begin{align}
\sum_{i=1}^m  (\la_{m+1} - \la_i)^p \leq  \dfrac{p \beta}{\gamma
N} \sum_{i=1}^m (\la_{m+1} - \la_i)^{p-1}  \la_i. \label{eq:3.20}
\end{align}
\end{cor}

\noindent Corollaries \ref{cor:3.6} and \ref{cor:3.7} follow from the facts about
%Corollaries 9 and 10
$\rho_i$ and $\Lambda_i$ given in conjunction with Corollaries \ref{cor:3.3}
and \ref{cor:3.4}, respectively.
%Corollaries 4 and 5

\section{The Case of a Schr\"{o}dinger-like Operator}  \label{schrodinger}
%Section IV
In this section, we consider an operator $H = A + V$ defined on $\mathcal{D}
\subset \mathcal{H}$, where $A$ and $V$ are self-adjoint operators,
$A= - \sum_{k=1}^N T_k^2$, and the $T_k$'s are skew-symmetric with domains
$T_k (\mathcal{D})$ satisfying $\mathcal{D} \equiv \mathcal{D}(A) \subset
\mathcal{D}(T_k)$ and $T_k (\mathcal{D}) \subset \mathcal{D} (T_k)$.  This
operator is modeled on the Schr\"{o}dinger operator.  We assume that the
spectrum of $H$ is discrete consisting of eigenvalues $\la_1 \leq \la_2 \leq \cdots$,
and we let $\{u_i\}_{i=1}^{\infty}$ be a complete orthonormal basis of eigenvectors
corresponding to $\{\la_i\}_{i=1}^{\infty}.$  We further take a family of symmetric
operators $\{B_k : H(\mathcal{D}) \rightarrow \mathcal{H} \}_{k=1}^{N}$ which
leave $\mathcal{D}$ invariant, such that $[T_\ell, B_k] u_i = \delta_{\ell k} u_i$.
As in Section \ref{sec:setup}, the quantities $\rho_i$ and $\Lambda_i$ are given
by
\[\rho_i = \sum_{k=1}^{N} \langle [H,B_k] u_i, B_k u_i\rangle,\]
and
\[
\Lambda_i = \sum_{k=1}^{N} \| [H,B_k] u_i \|^2.
\]
In obvious notation, we have $\rho_i = \rho_i^A + \rho_i^V$, corresponding to the
decomposition $H = A + V.$  The following theorem generalizes Theorem 4.1 of
\cite{AH2}.
\begin{thm}
\label{HS-Schro}  Suppose $[V, B_k] = 0$ for $1\leq k \leq N.$  Then $\rho_i = N$,
$\Lambda_i = 4 (\la_i - \langle V u_i, u_i \rangle).$  Moreover, for $p\le 2$
%for $0 \le p \le 2$
\begin{align}
\label{HS-Schro-eq1} \sum_{i=1}^m (\la_{m+1} - \la_i)^p \leq
\dfrac{4}{N} \sum_{i=1}^m (\la_{m+1} - \la_i)^{p-1} \big(\la_i -
\langle V u_i, u_i \rangle\big)
\end{align}
and for $p\ge 2$
\begin{align}
\sum_{i=1}^m (\la_{m+1} - \la_i)^p \leq \dfrac{2p}{N} \sum_{i=1}^m
(\la_{m+1} - \la_i)^{p-1} \big(\la_i - \langle V u_i, u_i
\rangle\big).
\label{HS-Schro-eq2}
\end{align}
\end{thm}
\smallskip \noindent {\em Proof.}
For the details of the calculations of $\rho_i$ and $\Lambda_i$,
see \cite{AH2}.  The rest follows via our previous considerations.
\hfill $\Box$

\begin{cor}
\label{cor:HS-Schro}
Suppose $V \ge M>0$.  Then, the inequalities in the previous theorem
reduce to
\begin{align}
\sum_{i=1}^m (\la_{m+1} - \la_i)^p \leq \dfrac{4}{N} \sum_{i=1}^m
(\la_{m+1} - \la_i)^{p-1} \big(\la_i -M\big) \text{ for } p\le 2
%\text{ for } $0 \le p \le 2
\label{HS-Schro-eq3}
\end{align}
and
\begin{align}
\sum_{i=1}^m (\la_{m+1} - \la_i)^p \leq \dfrac{2p}{N} \sum_{i=1}^m
(\la_{m+1} - \la_i)^{p-1} \big(\la_i - M\big) \text{ for } p\ge
2. \label{HS-Schro-eq4}
\end{align}
\end{cor}

\section{Comparing the Bounds} \label{sec:bounds}
%Section V
This section deals with the different bounds for $\la_{m+1}$
arising from the Harrell and Stubbe considerations \cite{HS} and
their extensions as detailed above.  We will assume that the
operators $A$ and $B_k$, $1\leq k \leq N,$ satisfy the conditions
(\ref{eq:3.70}) and (\ref{eq:3.80}) of Corollary \ref{cor:3.4},
namely
\[
\gamma \leq  [B_k, [A, B_k]]
\]
and
\[
- \sum_{k=1}^N [A, B_k]^2 \leq \beta A
\]
for some $\beta, \gamma >0$.

\subsection{Case  $p\le2$}
%\subsection{Case $0 \le p \le 2$}
%Should this say $0 \le p \le 2$?

We first treat the case $p \le2$, namely inequality (\ref{eq:3.91}) (or (\ref{eq:intro.6})
%Should this say $0 \le p \le 2$?
in the Introduction).  We assume $m\ge 2.$  For $m=1$ all bounds reduce to
\[\la_2 \le \left(1+ \dfrac{2 \beta}{\gamma N}\right) \ \la_1.\]
We set
\[f_p(\sigma)= \frac{1}{m} \ \sum_{i=1}^m (\sigma - \la_i)^p - \dfrac{2
\beta}{\gamma N} \ \frac{1}{m} \  \sum_{i=1}^m (\sigma -
\la_i)^{p-1} \ \la_i,\] for $\sigma \ge \la_m$. The {\em unique}
zero of $f_p(\sigma)$ larger than $\la_m$ is denoted by $\sigma_p$ (the existence
and uniqueness of $\sigma_p$ are addressed in Proposition \ref{prop:1} below).  It
%Proposition 13
can be thought of as a function of the moments $S_\ell$, for $\ell=1,
2, \cdots$, in the eigenvalues,
\begin{equation}
S_\ell = \dfrac{1}{m} \ \sum_{i=1}^m \la_i^{\ell}. \label{s1}
\end{equation}
This point of view becomes clear upon expansion in infinite series in
$\sigma$ (unless $p=1$ or $2$).  We first establish the following.
\begin{prop}
\label{prop:1}
For $0 \le p \le 2,$ there exists a unique root $\sigma_p$ of $f_p (\sigma)=0$ larger
than $\la_m$.
\end{prop}

\smallskip \noindent {\em Remark.}  The $p=2$ case of this proposition
is treated in detail in \cite{AH2} (see also \cite{A3}).  This proposition also holds
when $m=1,$ and then $\sigma_p = (1+\frac{2\beta}{\gamma N})\, \la_1$ (for all
$p$!).

\smallskip \noindent {\em Proof.}
{\em Existence.}  For $2\ge p>1$, $f_p(\la_m) \le 0$ by Corollary 5
with $m$ replaced by $m-1.$  This also holds, but as a strict inequality, for $p=1$,
since the inequality follows from the $p=1$ case of Corollary 5 (with $m$
replaced by $m-1$), but $f_p(\la_m)$ has one extra strictly negative term coming
from the $i=m$ term in the second sum in the definition of $f_p.$  Finally, for
$1>p\ge 0$,
\[\lim_{\sigma \rightarrow \la_m^{+}} f_p(\sigma) = -\infty.\]
Moreover, for $p > 0$ \[\lim_{\sigma \rightarrow \infty} f_p(\sigma) = \infty,\]
while for $p=0$
\[\lim_{\sigma \rightarrow \infty} f_0(\sigma) = 1.\]

Hence, by continuity, the existence of a zero larger than
$\la_m$ of $f_p(\sigma)=0$ (where $0 \le p \le 2$) is guaranteed.  In fact we
can say more.

{\em Uniqueness.}  First, observe that \[f_0(\sigma) = 1 - \frac{2
\beta}{\gamma N} \ \frac{1}{m} \  \sum_{i=1}^m \frac{\la_i}{\sigma
- \la_i}.\]
Hence $f_0'(\sigma)>0$ and $f_0(\sigma)$ is
monotonically increasing from $-\infty$ to $1.$  This establishes
the uniqueness of $\sigma_0$ (this is the Hile-Protter
bound for $\la_{m+1}$ derived in \cite{HP} (see also \cite{A3},
\cite{A2}, \cite{AH2}, \cite{HM2}, \cite{HM1}, \cite{HS},
\cite{Michel}, \cite{Pr4}.  For $0<p\le 1$,
\[f'_p(\sigma) = \frac{1}{m} \ \Big( p \ \sum_{i=1}^m \ (\sigma -
\la_i)^{p-1}- (p-1) \frac{2 \beta}{\gamma N} \ \sum_{i=1}^m
(\sigma - \la_i)^{p-2} \ \la_i\Big).\]
(The second term is identically $0$ if $p=1$.)  Since $p-1 \le 0$ for $0 < p \le 1, $
$f'_p(\sigma)> 0$ for $\sigma>\la_m$ and the uniqueness of $\sigma_p$ is
established in this case.  Note that this handles the $p=1$ case since, as already
observed, $f_1(\la_m) < 0$ (strict inequality).  This case can also be treated via
explicit and elementary calculation (as can the $p=2$ case).

For $1<p< 2,$ we note that $f'_p(\sigma)$ is not clearly of one sign as before,
since $f'_p(\sigma) \rightarrow -\infty$ as $\sigma \rightarrow \la_m^+$, while
$f'_p(\sigma) \rightarrow \infty$ as $\sigma \rightarrow \infty$.  We therefore have
recourse to a convexity argument.  Differentiating, it becomes clear that in this case
$f''_p(\sigma)> 0$ for $\sigma>\la_m$.  Hence $f'_p(\sigma)$ is strictly increasing
from $-\infty$ (value at $\la_m$ for $1<p<2$) to $\infty$ (value at $\infty$).  One
can then find a unique $\xi_p>\la_m$ for which $f'_p(\xi_p)=0$.  Moreover,
$f'_p(\sigma) < 0$ for $\la_m<\sigma < \xi_p$ and $f'_p(\sigma) > 0$ for
$\sigma > \xi_p.$  The uniqueness of $\sigma_p$ is therefore ascertained with
$\la_m<\xi_p<\sigma_p$.  Finally, the case of $p=2$ follows easily in much the
same way as for $1 < p < 2$ using now the fact that $f_2(\sigma)$ is a quadratic
in $\sigma$ with second order term $\sigma^2.$  Thus $f_p$ is again concave up
and the result follows.
\hfill $\Box$

Since $f_p(\la_{m+1}) \le 0$ (viz., (\ref{eq:3.91})), it obtains that
$\la_{m+1} \le \sigma_p$ for any $0\le p\le 2.$  Moreover, since
$f_p((1+ \frac{2 \beta}{\gamma N}) \la_m) > 0$, we have the
inequalities

\begin{align}
\label{eq:quadrature}
\la_m \le \la_{m+1} \le \sigma_p < \left(1+ \frac{2 \beta}{\gamma N} \right) \la_m.
\end{align}

We are now ready to prove the statement announced in \cite{AH2}:
{\em ``$\sigma_p$ improves with $p$, for $p\le 2.$''} This is
contained in the following theorem.

\begin{thm} \label{thm:4.1}
\begin{align}
\la_{m+1} \le \sigma_{p_2} \le \sigma_{p_1} \qquad {\text if}
\qquad 0\le p_1 \le p_2 \le 2.
\label{eq:4.1}
\end{align}
\end{thm}

\smallskip \noindent {\em Proof.}  This is done in several reductions.  We
observe that the statement $\sigma_{p_2} \le \sigma_{p_1}$ is equivalent
to showing that $f_{p_1} (\sigma_{p_2}) \le 0$ (since $f_{p_1} (\sigma)$ is
below the $\sigma$-axis on $(\la_m, \sigma_{p_1})$).  That is
\[
\sum_{i=1}^m (\sigma_{p_2} - \la_i)^{p_1}
\le \dfrac{2 \beta}{\gamma N} \sum_{i=1}^m (\sigma_{p_2} - \la_i)^{p_1-1}
\ \la_i.\]

Since
\[\dfrac{2 \beta}{\gamma N} = \dfrac{\sum_{i=1}^m (\sigma_{p_2} - \la_i)^{p_2}}
{\sum_{i=1}^m (\sigma_{p_2} - \la_i)^{p_2-1} \ \la_i},\]
the statement of the theorem is then equivalent to

\begin{align}
\dfrac{\sum_{i=1}^m (\sigma_{p_2} - \la_i)^{p_1}}{\sum_{i=1}^m (\sigma_{p_2} - \la_i)^{p_1-1} \ \la_i}
\le \dfrac{\sum_{i=1}^m (\sigma_{p_2} - \la_i)^{p_2}}
{\sum_{i=1}^m (\sigma_{p_2} - \la_i)^{p_2-1} \ \la_i}.
\label{eq:4.2}
\end{align}
Or
\begin{align}
\sum_{i=1}^m (\sigma_{p_2} - \la_i)^{p_1} \ \sum_{i=1}^m (\sigma_{p_2} - \la_i)^{p_2-1} \ \la_i
\le \sum_{i=1}^m (\sigma_{p_2} - \la_i)^{p_2} \ \sum_{i=1}^m (\sigma_{p_2} - \la_i)^{p_1-1} \ \la_i.
\label{eq:4.3}
\end{align}

We now use the following version of the ``Chebyshev Inequality''
(see, for example, p.\ 43 of \cite{HLP}).
\begin{lem} [{\bf Weighted Reverse Chebyshev Inequality}]
\label{lemma:cheb}  Let $\{a_i\}_{i=1}^m$ and $\{b_i\}_{i=1}^m$ be
two oppositely ordered real sequences, and let $\{w_i\}_{i=1}^m$
be a sequence of nonnegative weights.  Then the following
inequality holds
\begin{align}
\sum_{i=1}^m w_i \ \sum_{i=1}^m w_i \ a_i b_i \le \sum_{i=1}^m w_i
\ a_i \ \sum_{i=1}^m w_i \ b_i.
\label{eq:4.4}
\end{align}
\end{lem}

Inequality (\ref{eq:4.3}) is then a corollary to this lemma with
$w_i= (\sigma_{p_2} - \la_i)^{p_1}$, $a_i =
\frac{\la_i}{\sigma_{p_2} - \la_i}$, and $b_i=(\sigma_{p_2} -
\la_i)^{p_2-p_1}.$  The sequence $\{a_i\}$ is increasing, while the
sequence $\{b_i\}$ is decreasing by virtue of the fact that $p_2
\ge p_1.$  Hence the result of the theorem.
\hfill $\Box$

\noindent {\em Remarks.}

1. This theorem contains the statement announced by H.\ C.\ Yang \cite{Y} that
his inequality ($p=2;$ also referred to as ``Yang 1'') implies an ``averaged''
version of this inequality ($p=1;$ also referred to as ``Yang 2'') which in turn
implies the Hile-Protter result ($p=0$).  In \cite{A2}, this statement is
summarized in the implication that (for each $m=1, 2, \dots$)

\[\text{ Yang } \, 1 \implies \text{ Yang } \, 2   \implies \text{ Hile-Protter }.\]
%Check that \text is doing the right thing here.  I was trying to get away from the
%first letters in the {\text     }s not being in italics while everything else was, plus
%the extra spaces that were occurring around the hyphen in ``Hile-Protter''.

2. A proof of this result is not given in \cite{Y}.  Proofs are given
in \cite{A3} and \cite{AH2} (see also \cite{A2}).
%In this paper, we presented a proof which partially emulates that
%of \cite{AH2}.
Our proof here is basically that of \cite{AH2}, but gives a more
general result.  This theorem shows that of the class of
Harrell-Stubbe-type inequalities with $p \le 2$,
%$0 \le p \le 2$,
the optimum obtains when $p=2$ (H.\ C.\ Yang).

3. The PPW inequality is of course weaker than the HP inequality.
%It provides a tighter bound for $\sigma_p$ than (\ref{eq:quadrature}), viz.,
%What does ``It'' refer to here?  We should clean up this sentence to make this
%reference (i.e., what the pronoun ``It'' refers to) clearer.
Thus the HP inequality provides a tighter bound for $\sigma_p$ than the bound
\begin{align}
\label{eq:ppw}
\sigma_p \le \la_m + \frac{2 \beta}{\gamma N} \ S_1,
\end{align}
which is the PPW inequality in this setting.  Also, it is perhaps worth noting that
the HP inequality provides a tighter bound for $\sigma_p$ than that given in
\eqref{eq:quadrature}, viz.,
%The reference here (for \eqref{eq:quadrature} does not come through.  We need
%to check this out and correct whatever is wrong.
\[\sigma_p <  \left(1 + \dfrac{2 \beta}{\gamma N} \right) \, \la_m\]
(which is itself a simple consequence of the PPW inequality).
%(compare with \eqref{eq:quadrature}).

4. In terms of the ``moments'', $S_{\ell}= \frac{1}{m} \
\sum_{i=1}^m \la_i^{\ell}$, Yang 2 reads,
\begin{align}
\label{yang2}
\sigma_1 = \left(1+ \frac{2 \beta}{\gamma N}\right) S_1,
\end{align}
while Yang 1 translates as
\begin{align}
\label{yang1}
\sigma_2 = \left(1+ \frac{\beta}{\gamma N}\right) S_1 +\left\{ \left(1+ \frac{\beta}{\gamma N}\right)^2 S_1^2 -
\left(1+ \frac{2 \beta}{\gamma N}\right) S_2  \right\}^{1/2}.
\end{align}

\subsection{Case $p\ge2$}

We now turn our attention to Harrell and Stubbe's second extension
of H.\ C.\ Yang's result, namely Corollary \ref{cor:3.7}
%Corollary 10
(where
$p\ge 2$).  We will show that bounds for $\la_{m+1}$ provided by
(\ref{eq:3.20}) (or (\ref{eq:intro.5}) in the Introduction)
obtained when $p\ge 2$ are weaker than those for $p=2$.  In fact
they get worse monotonically with increasing $p$.

In this case, the function $f_p(\sigma)$ takes a slightly altered form,
which we denote by $\tilde{f}_p(\sigma)$:
\begin{align}
\tilde{f}_p(\sigma) = \frac{1}{m} \sum_{i=1}^m (\sigma - \la_i)^p -
\dfrac{p\beta}{\gamma N} \ \frac{1}{m}  \ \sum_{i=1}^m \
(\sigma - \la_i)^{p-1} \,\la_i.
\label{eq:4.6.2}
\end{align}
We will denote this function by $\tilde{f}_{p,m}(\sigma)$ in case the
explicit dependence of $\tilde{f}_p(\sigma)$ on $m\ge 1$ is required.
The existence of a root of $\tilde{f}_p(\sigma)=0$ greater than or equal
to $\la_m$ is guaranteed.  This is because of $\tilde{f}_{p,m}(\la_m) =
\frac{m-1}{m} \tilde{f}_{p,m-1}(\la_m) \le 0$ (viz., (\ref{eq:3.20})) and
$\lim_{\sigma \rightarrow \infty} \tilde{f}_p(\sigma) = \infty$, and because
$\tilde{f}_p(\sigma)$ is continuous on $(\la_m, \infty).$
%We will deal only with the case when this root is unique.  We will
%denote this unique root by $\tilde{\sigma}_p.$  Under this definition,
%it is clear that $f_p(\sigma) \le 0$ for $\la_m \le \sigma \le
%\tilde{\sigma}_p$.

We handle this case somewhat differently from how we handled the
case for $0 \le p \le 2.$  In that case we established the existence of a
unique root of $f_p(\sigma)$ greater than $\la_m.$  In the present case
we do not establish uniqueness (although it may well obtain) but rather
define (existence follows from our comments above) a root $\tilde{\sigma}_p$
of $\tilde{f}_p(\sigma)$ which is greater than or equal to $\la_m$ and serves
our purposes.  For $p \ge 2$ we define $\tilde{\sigma}_p$ via
\begin{align}
\tilde{\sigma}_p = sup \{s \ge \la_m | \tilde{f}_p(\sigma) \le 0 \text{ for } \la_m \le
\sigma \le s \}.
\end{align}
By continuity, $\tilde{\sigma}_p$ is in fact realized as a maximum over
the set given on the right, and we have $\tilde{f}_p(\sigma) \le 0$ for
$\la_m \le \sigma \le \tilde{\sigma}_p.$ Indeed, it must also be true that
$\tilde{f}_p(\tilde{\sigma}_p)=0,$ i.e., that $\tilde{\sigma}_p$ is a root of
$\tilde{f}_p(\sigma)$ which is greater than or equal to $\la_m.$  Note, too,
that from the fact that $\tilde{f}_2(\sigma)$ is a quadratic with leading term
$\sigma^2$ and that $\tilde{f}_2(\sigma) \le 0$ on $[\la_m,\la_{m+1}],$ it is
clear that $\tilde{\sigma}_2$ as defined above is identical to $\sigma_2$ as
defined previously (for the case of $0 \le p \le 2),$ which is also just the
explicit upper bound for $\la_{m+1}$ coming from Yang's first (or main)
inequality, that is, the expression given on the right-hand side of
(\ref{yang1}).

As remarked earlier in a similar context, $\tilde{\sigma}_p$ can be thought
of as a function of the moments in the first $m$ eigenvalues
providing an upper bound for $\la_{m+1}.$  In this case,
$\tilde{f}_p(\sigma)$ takes the form
\begin{align}
\tilde{f}_p(\sigma) = \sigma^p + \sum_{k=1}^{N_p} (-1)^k \binom{p}{k}
\left(1+ \frac{\beta k}{\gamma N} \right) S_k \sigma^{p-k}
\notag
\end{align}
where $S_k$ is defined by (\ref{s1}) and $N_p=p$ if p is an integer
and $N_p=\infty$ if $p$ is not an integer.  By convention, the binomial
coefficient $\binom{p}{k}$ denotes \[\frac{p(p-1)(p-2) \cdots (p-k+1)}{k!}\]
even when $p$ is not an integer.

To proceed, we need to know that $\tilde{\sigma}_p$ as defined
above for $p > 2$ really does provide a bound for $\la_{m+1}$ (for
this it is not enough to know Corollary \ref{cor:3.7}, i.e.,
ineq.\ \eqref{eq:3.20},
%Corollary 10
%ineq.\ (48)
since our definition of $\tilde{\sigma}_p$
does not preclude the possibility that $\tilde{f}(\sigma) = 0$ has further roots
beyond $\tilde{\sigma}_p,$ or further places where $\tilde{f}_p(\sigma) < 0,$
and that $\la_{m+1}$ is then somewhere to the right of $\tilde{\sigma}_p$ as
defined above).

There are (at least) three ways we could think to proceed at this point.

%The itemization below needs some tex commands for it to come out
%the way I want it, I suppose.
(1)  Use the $p=2$ case and a technique of Aizenman and Lieb
\cite{Aizen-Lieb} (see also \cite{Hund-Sim}, \cite{Lap}) to show that
$[\la_m,\tilde{\sigma}_2] \subset [\la_m, \tilde{\sigma}_p]$ for $p > 2$ and
hence that $\la_{m+1} \le \tilde{\sigma}_2 \le \tilde{\sigma}_p$ for
$p \ge 2,$ by our definition of $\tilde{\sigma}_p.$

(2)  Specialize to the case of ineq.\ \eqref{eq:intro.5} from our Introduction, i.e., to the case
%ineq.\ (2)
of the Laplacian (and certain generalizations), where Harrell and Stubbe
\cite{HS} have already provided results implying that $\la_{m+1} \le \tilde{\sigma}_p$
(and indeed that $\tilde{f}_p(\sigma) \le 0$ on $[\la_m, \la_{m+1}]).$  For these
results we refer to ineq.\ (14) in Theorem 9 on p.\ 1805 and, in particular, the
conditions that go with it.  Specifically, their results show that (for $\tilde{f}_p(\sigma)$
as in \eqref{eq:4.6.2} but with $\beta p/\gamma N$ replaced by $2p/n)$ $\tilde{\sigma}_p \ge
%eq.\ (61)
\la_m + (p/2)(\la_{m+1} - \la_m)$ and since $\la_m + (p/2)(\la_{m+1}-\la_m) > \la_{m+1}$ for
$p > 2,$ the desired result follows.  One can consult \cite{HS}, \cite{AH4} for generalizations of
%I need to supply the ref. here.  What's the appropriate ref. for this?  Is \cite{AH4} an appropriate
%reference to one of our papers?  Drop the reference to \cite{AH4} if it's not appropriate.
$-\Delta$ to which the Harrell-Stubbe results are already known to apply.

(3)  Extend the approach and methods of Harrell and Stubbe \cite{HS} so that we
know that, for the operators considered here and for $\tilde{f}_p(\sigma)$ as defined by
\eqref{eq:4.6.2}, $\tilde{f}_p(\sigma) \le 0$ for all $\sigma \in [\la_m, \la_{m+1}]$ and $p \ge 2.$
%eq.\ (61)
Thus, under this approach one would beginb by seeking a version of Harrell and Stubbe's
Theorem 9 (p.\ 1805 of their paper) that applies in our general setting.

In what follows we will follow (1) since it gives the most self-contained approach from
our chosen point of view.  One could also build on (2), which puts one farther along with
the problem at the start, but, as mentioned above, leads to a more restricted result.
Finally, (3) would probably also work, and lead to results analogous to and as general as
those of (1) (even, perhaps, to results which are a bit stronger), but as we have not worked
through the details of this we leave it aside.

\begin{thm} \label{thm:4.2}  Suppose $p \ge 2.$  Then $\tilde{f}_p(\sigma) \le 0$ for all
$\sigma \in [\la_m, \tilde{\sigma}_2]$ and hence, since $\la_{m+1} \le \tilde{\sigma}_2,
\la_{m+1} \le \tilde{\sigma}_p.$  Moreover $\tilde{\sigma}_p \ge \tilde{\sigma}_2
(=\sigma_2).$
\end{thm}

\smallskip \noindent {\em Proof.}  We know that
\begin{align}
\label{eq:4.2.1}
%eq.\ (62)  (It's possible that this number (i.e., 62) may change.)
m\tilde{f}_2(\sigma) =   \sum_{i=1}^{m} (\sigma-\la_i)^2 -
\frac{2 \beta}{\gamma N} \sum_{i=1}^{m}(\sigma-\la_i)\,\la_i \, \le \, 0
\end{align}
for all $\sigma \in [\la_m,\tilde{\sigma}_2],$ and, in particular, for all
$\sigma \in [\la_m,\la_{m+1}].$  Because this holds for all values of $m$ (and specifically
for $1,2, \dots, m$ replacing $m$ above), if we introduce the notation
\begin{align}
\left(\la - t\right)_{+}=
\begin{cases} \la-t & \text{ if } t \le \la, \cr
0 & \text{ if } t>\la, \cr
\end{cases}
\end{align}
%Another rendering of this could be as follows:
%\[\nonumber
%(\la-t)_{+}=\left\{
%\begin{array}{c}
%\displaystyle{\la-t,\;\;{\mbox{if}} \;\; t \le \la,}\\
%\displaystyle{0 \ ,\;\;{\mbox{if}} \;\;\;t\ge \la.}
%\end{array}
%\right.
%\]
then \eqref{eq:4.2.1} extends to all $\sigma \le \tilde{\sigma}_p$ as
%eq.\ (62)  (It's possible that this number (i.e., 62) may change.)
\begin{align}
\label{eq:4.2.2}
%ineq. (63)  (It's possible that this number (i.e., 63) may change.)
\sum_{i=1}^m\,(\sigma-\la_i)^2_+ - \dfrac{2 \beta}{\gamma N} \sum_{i=1}^m\,(\sigma-\la_i)_+\, \la_i
\le 0.
\end{align}
Note that each time $\sigma$ passes below a $\la_i$ another term drops away (on both sides),
leaving us with a variant of ineq.\ \eqref{eq:4.2.1} where the only change is that $m$ is less.  And finally,
%ineq.\ (62)  (It's possible that this number (i.e., 62) may change.)
when $\sigma$ crosses $\la_1$ we are left with the trivial inequality $0 \le 0.$

We now rewrite \eqref{eq:4.2.2} as
%eq.\ (63)  (It's possible that this number (i.e., 63) may change.)
\begin{align}
\label{eq:4.2.3}
%ineq. (64)  (It's possible that this number (i.e., 64) may change.)
\sum_{i=1}^m\,(\sigma-\la_i-r)^2_+ - \dfrac{2 \beta}{\gamma N} \sum_{i=1}^m\,(\sigma-\la_i-r)_+\, \la_i
\le 0,
\end{align}
which holds for all $r \ge 0$ if $\sigma \le \tilde{\sigma}_2.$  In particular we consider $\sigma \in
(\la_1, \tilde{\sigma}_2].$  If we integrate this inequality against $r^{p-3}$ for $0<r<\infty$ we can
use the beta function integral
\begin{align}
\label{eq:4.2.4}
%eq.\ (65)  (It's possible that this number (i.e., 65) may change.)
%Perhaps this eq. doesn't need to be labeled?
B(s,t)= \int_0^1\,u^{s-1}\,(1-u)^{t-1}\,du   \text{ for }  s, t >0
&=\dfrac{\Gamma(s)\,\Gamma(t)}{\Gamma(s+t)}
\end{align}
to evaluate the integrals (this is the ``trick'' of Aizenman and
Lieb \cite{Aizen-Lieb}). We have (for $\alpha > -1, p>2$
)
\begin{align} \label{eq:4.2.5}
%eq.\ (66)  (It's possible that this number (i.e., 66) may change.)
%Perhaps this eq. doesn't need to be labeled?
\int_0^\infty(\sigma-\la_i-r)^\alpha_+ \, r^{p-3}\, dr=0 \quad
\text{ if } \sigma-\la_i \le 0. \notag
\end{align}
This integral reduces to
\begin{align}
\int_0^{\sigma-\la_i} ((\sigma-\la_i)-r)^\alpha \, r^{p-3}\, dr
\quad \text{ if } \sigma-\la_i > 0.
\end{align}
Changing variables via $r=(\sigma-\la_i)u$ for $\sigma-\la_i$ positive, we arrive at
\begin{align}
\label{eq:4.2.6}
%eq.\ (67)  (It's possible that this number (i.e., 67) may change.)
%Perhaps this eq. doesn't need to be labeled?
\int_0^\infty \, (\sigma-\la_i-r)^\alpha_+ \, r^{p-3} \, dr = (\sigma -\la_i)^{\alpha+p-2}_+ \, B(p-2,\alpha+1)
\end{align}
and hence ineq.\ \eqref{eq:4.2.3} becomes (for $p > 2$)
%ineq.\ (64)  (It's possible that this number (i.e., 64) may change.)
%Ineq. (64) is the ineq. that appears after ``We now rewrite (63) as . . .'' above.
\begin{align}
\label{eq:4.2.7}
%eq.\ (68)  (It's possible that this number (i.e., 68) may change.)
%Perhaps this eq. doesn't need to be labeled?
B(p-2,3) \, \sum_{i=1}^m \, (\sigma-\la_i)^p_+ \, -\dfrac{2\beta}{\gamma N} \, B(p-2,2) \, \sum_{i=1}^m \,
(\sigma - \la_i)^{p-1}_+ \, \la_i \le 0,
\end{align}
or, since
$$\dfrac{B(p-2,2)}{B(p-2,3)}=\dfrac{\Gamma(p-2)\Gamma(2)}{\Gamma(p)} \,
\dfrac{\Gamma(p+1)}{\Gamma(p-2)\Gamma(3)}$$
and $\Gamma(p+1)=p \, \Gamma(p),$ $\Gamma(1)=1,$
\begin{align}
\label{eq:4.2.8}
%ieq.\ (69)  (It's possible that this number (i.e., 69) may change.)
%This ineq. needs to be labeled appropriately.  For now I'm calling it ineq. (69) (formerly (65)) in my text.
\sum_{i=1}^m \, (\sigma-\la_i)^p_+ \, -\dfrac{p\beta}{\gamma N} \, \sum_{i=1}^m \,
(\sigma - \la_i)^{p-1}_+ \, \la_i \le 0,
\end{align}
for all $\sigma \in (\la_1,\tilde{\sigma}_2].$

Thus the inequality $\tilde{f}_p(\sigma) \le 0$ holds for all $\sigma \in [\la_m, \tilde{\sigma}_2]$
(and similarly when $m$ is replaced by any positive integer if $\tilde{\sigma}_2$ is understood
as $\tilde{\sigma}_{2,m},$ the root $\tilde{\sigma}_2$ when $\tilde{f}_2(\sigma) = \tilde{f}_{2,m}(\sigma);$
in particular, we have $\tilde{f}_{2,k} \le 0$ for $\sigma \in [\la_k,\la_{k+1}]$ since
$[\la_k,\la_{k+1}] \subset [\la_k, \tilde{\sigma}_{2,k}]$).  The definition of $\tilde{\sigma}_{p,m}$
now implies that $\la_{m+1} \le \tilde{\sigma}_{2,m} \le \tilde{\sigma}_{p,m}$ (and, in fact, that
$\tilde{f}_{p,m}(\sigma) \le 0$ on $[\la_m, \tilde{\sigma}_{p,m}] \supset
%Is \supset right here?  I intend this as a set containment, but in the sense opposite to \subset.
[\la_m,\tilde{\sigma}_{2,m}]$), or dropping again the $m$ subscript on $\tilde{\sigma}_p,$
$\la_{m+1} \le \tilde{\sigma}_2 \le \tilde{\sigma}_p$ for $p \ge 2,$ which is the final conclusion we
wished to draw.
\hfill $\Box$
%I'd like to put an ``end of proof'' box here.  I guess this should do it (\hfill $\Box$).

Thus the $p=2$ bound for $\la_{m+1}$ equals or surpasses all the bounds $\tilde{\sigma}_p$ for
$p \ge 2$ coming from ineq.\ \eqref{eq:3.20}
%ineq.\ (48)  (It's possible that this number (i.e., 48) may change.)
via our definition of the $\tilde{\sigma}_p$'s.  This is certainly enough, from one point of view, to
dismiss the inequality \eqref{eq:3.20}
%ineq.\ (48)  (It's possible that this number (i.e., 48) may change.)
for all $p \ge 2$ from further consideration but we cannot resist drawing one final conclusion from
the Aizenman-Lieb technique.
\begin{thm}
\label{thm:4.3}
%I'll call this Theorem 17.  Probably the boldfacing I'm doing here will need to come out.
%{\bf Theorem 17.}
Suppose $q \ge p \ge 2.$  Then
\begin{align}
\label{eq:4.3.1}
%ineq.\ (70)  (It's possible that this number (i.e., 70) may change.)
\la_{m+1} \le \tilde{\sigma}_2 \le \tilde{\sigma}_p \le \tilde{\sigma}_q.
\end{align}
\end{thm}
%Probably my boldfacing of this ``Proof'', etc. below also needs adjustment.
%{\bf Proof.}
\smallskip \noindent {\em Proof.}  One proceeds from ineq.\ \eqref{eq:4.2.8}
%ineq.\ (69) (formerly (65))
much as we did from ineq.\ \eqref{eq:4.2.2} above, first putting it in the form
%ineq.\ (63)  (It's possible that this number (i.e., 63) may change.)
\begin{align}
\label{eq:4.3.2}
%ineq.\ (71)  (It's possible that this number (i.e., 71) may change.)
%Perhaps this ineq. doesn't need to be labeled?
\sum_{i=1}^m \, (\sigma-\la_i-r)^p_+ \, - \, \dfrac{p\beta}{\gamma N} \,
\sum_{i=1}^m \, (\sigma-\la_i-r)^{p-1}_+ \, \la_i \, \le \, 0
\end{align}
which we know to hold for all $\sigma \le \tilde{\sigma}_p$ and $r \ge 0.$  One then integrates
in $r$ much as before, except that this time one multiplies by $r^{q-p-1}$ (for $q > p$) before
integrating from 0 to $\infty.$  This leads to
\begin{align}
\label{eq:4.3.3}
%ineq.\ (72)  (It's possible that this number (i.e., 72) may change.)
%Perhaps this ineq. doesn't need to be labeled?
\sum_{i=1}^m \, (\sigma-\la_i)^q_+ \, - \, \dfrac{p\beta}{\gamma N} \, \dfrac{B(q-p,p)}{B(q-p, p+1)} \,
\sum_{i=1}^m \, (\sigma-\la_i)^{q-1}_+ \, \la_i \, \le \, 0
\end{align}
for all $\sigma \le \tilde{\sigma}_p,$ which we can extend to all $\sigma \le \tilde{\sigma}_q$ by
how we defined the $\tilde{\sigma}_p$'s.  Noting that
$$\dfrac{B(q-p,p)}{B(q-p,p+1)}=\dfrac{\Gamma(p) \, \Gamma(q+1)}{\Gamma(p+1) \, \Gamma(q)}
=\dfrac{q}{p}$$
we see that we have arrived at
\begin{align}
\label{eq:4.3.4}
%ineq.\ (73)  (It's possible that this number (i.e., 73) may change.)
%Perhaps this ineq. doesn't need to be labeled?
\sum_{i=1}^m \, (\sigma-\la_i)^q_+ \, - \, \dfrac{q\beta}{\gamma N} \,
\sum_{i=1}^m \, (\sigma-\la_i)^{q-1}_+ \, \la_i \, \le \, 0,
\end{align}
which is what we sought, since we have that $\tilde{f}_q(\sigma) \le 0$ for
$\sigma \in [\la_m, \tilde{\sigma}_p]$ and, extending via the definition of $\tilde{\sigma}_q,$ for all
$\sigma \in [\la_m,\tilde{\sigma}_q].$

Thus we have $\tilde{\sigma}_2 \le \tilde{\sigma}_p \le \tilde{\sigma}_q$ for $q \ge p$ and since
we know that $\la_{m+1} \le \tilde{\sigma}_2$ this completes the proof of the theorem.
\hfill $\Box$
%I'd like to put an ``end of proof'' box here.  I guess this should do it (\hfill $\Box$).

%We now delete up to the end of the old proof of Theorem 17.
%Then we take up again with the old text at the middle of p. 23 (i.e., with ``We end this section by . . .''),
%just after the ``end of proof'' box.
%We could leave in the following paragraph, which would then follow the end of the proof just above.
A small remark here is that there is a nice identity $\tilde{f}'_p(\sigma)=p \, \tilde{f}_{p-1}(\sigma),$
showing that zeros of $\tilde{f}_{p-1}$ are critical points of $\tilde{}_p.$  While this allows one to start
analyzing the behavior of $\tilde{f}_p$ based upon that of $\tilde{f}_{p-1}$ we were not able to build
a general approach along these lines.  And, at best, even if successful this approach would only allow
comparisons of $\tilde{\sigma}_p$'s for values of $p$ differing by an integer.

We end this section by mentioning that the (standard) Reverse
Chebyshev Inequality implies (via an argument similar to that used
in our proof of Theorem \ref{thm:4.1} above) that
%Theorem 14
$g_m (\tilde{\sigma}_p) \geq
\dfrac{m \gamma N}{p \beta}$ for
\[g_m(\sigma) = \sum_{i=1}^m \dfrac{\la_i}{\sigma- \la_i}.\]
This holds for any choice of $\tilde{\sigma}_p$ for $p \ge 2.$
Thus, we have the upper estimate
\[\tilde{\sigma}_p \leq \la_m + \dfrac{p \beta}{\gamma N} S_1\]
by replacing the quantities $\tilde{\sigma}_p -\la_i$ by the
smallest, i.e., $\tilde{\sigma}_p -\la_m$ in the expression of
$g_m(\tilde{\sigma}_p)$ (note that we already have
$\tilde{\sigma}_p \ge \la_m$).  This bound is in the spirit of
the PPW bound (\ref{eq:intro.1}) (cf.\ also \eqref{eq:ppw})
except for a $p$ in place of a 2 on the right-hand side.

It is not clear at this stage whether the Harrell-Stubbe inequality
is stronger than that of Hile-Protter (the $p=0$ case of \eqref{eq:3.91}
in its generalized form, or
%ineq.\ (30) (see Corollary 5)
\eqref{eq:intro.hp} originally) for all $p > 2$ or not.  It surely is,
by continuity, for some range of $p$'s just larger than 2.

\section{Applications to Physical and Geometric Problems}
\label{applications}

In this section, we illustrate some applications of the abstract
formulation described earlier.  Physical and geometric problems are
considered.  Our results improve earlier bounds for various
eigenvalue problems by Harrell and Michel \cite{HM2}, for
eigenvalues of domains in $\sz^2$ and $\hz^2$, as well as other
bounds by Hook \cite{Ho2}.  The general strategy, as explained in
\cite{AH4} (see also \cite{A3}, \cite{Ha}, \cite{Ha2}, \cite{HM2},
\cite{HM1}, \cite{Ho2}, \cite{LP}, \cite{Lee}, \cite{Leung},
\cite{Li}, \cite{Y}, \cite{YY}), is to write the operator $A$ in
the form $A=- \sum_{k=1}^N T_k^2 + V$, where the $T_k$'s are
skew-symmetric.  The auxiliary symmetric operators $B_k$ are chosen
such that $[T_\ell, B_k]= \delta_{\ell k}$ and $[V, B_k] = 0$.  The
``potential'' $V$ is either zero or appropriately bounded below.
Sometimes it is more appropriate to reduce to a situation like
that in Corollary \ref{cor:3.4} of Section III.  Once this is done, a family
of new inequalities of the Harrell-Stubbe-type is obtained for the
eigenvalue problem at hand.  We illustrate this via several examples.

\subsection{Classical PPW, HP, and Yang Inequalities for the Fixed Membrane}
%Subsection A
\label{classical} For the classical ``fixed membrane'' problem
described in Section \ref{intro}, $A = -\D$, $T_j =
\dfrac{\partial}{\partial x_j}$ and $B_j = x_j$, for $1\leq j \leq
n$ are the appropriate choices. We have
\[A = - \sum_{j=1}^n T_j^2,\] and
\[[T_\ell, B_k] = \delta_{\ell k}.\] Under the Dirichlet boundary
conditions of the problem, the $T_j$'s are skew-symmetric with
respect to the inner product
\[\langle u, v \rangle = \int_{\Om} u v \, dx.\]
The classical inequalities of PPW, HP, and Yang then follow straightforwardly
via the results presented in Sections III and IV, as do their Harrell-Stubbe-style
generalizations.

\subsection{The Inhomogeneous Membrane Problem} \label{inhomogeneous}
%Subsection B
This is of course a generalization of the  fixed membrane problem
in the previous section.  In this case, the density $q(x)$ of the
membrane is not uniform over $\Omega \subset \rz^n$. The
eigenvalue model for this problem is given by
\begin{align}
 - \Delta u & = \lambda \, q(x) \, u
\text { in } \Omega, \nonumber \\
u & = 0 \text { on } \partial \Omega. \label{eq:wdlp.1}
\end{align}
We assume $0<q_{min} \le q(x) \le q_{max} <\infty$.  The operator
$A$ takes the form $A = \dfrac{-\Delta}{q(x)}$.  It is symmetric
with respect to the inner product $\langle u , v \rangle_q =
\int_{\Om} u(x) v(x) q(x) dx$. The real eigenfunctions $\{ u_i
\}_{i=1}^{\infty}$ satisfy
\[\int_{\Om} u_i(x) u_j(x) q(x) dx = \delta_{ij},\]
and the eigenvalues $\{ \la_i \}_{i=1}^{\infty}$ are given by
$\la_i = \int_{\Om} |\Q u_i|^2 dx$. With $B_k= x_k$, for $1\le
k\le n$, one is led to $\rho_i\ge \frac{n}{q_{max}}$ and
$\Lambda_i\le \frac{4 \la_i}{q_{min}}$.  Hence, we have the
following extension of a result of Ashbaugh \cite{A3}.
\begin{thm} \label{inhomogeneous-thm}
The eigenvalues of the inhomogeneous membrane problem with
density function $0<q_{min} \le q(x) \le q_{max} <\infty$ satisfy the
inequalities
\begin{align}
\sum_{i=1}^m (\la_{m+1} - \la_i)^p \leq \dfrac{4}{n} \, \dfrac{
q_{max}}{q_{min}} \sum_{i=1}^m (\la_{m+1} - \la_i)^{p-1}  \la_i
\text{ for } p\le 2
%\text{ for } 0 \le p \le 2
\label{HS-inh1}
\end{align}
and
\begin{align}
\sum_{i=1}^m (\la_{m+1} - \la_i)^p \leq \dfrac{2 p}{n} \, \dfrac{
q_{max}}{q_{min}} \sum_{i=1}^m (\la_{m+1} - \la_i)^{p-1}  \la_i
\text{ for } p\ge 2. \label{HS-inh2}
\end{align}
\end{thm}
\noindent {\em Remark.}  Ashbaugh's result (see Section 4 of
\cite{A3}) is a refinement and strengthening of a result first
proved by Cheng \cite{Cheng} in the context of a minimal
hypersurface $\Omega$ in $\rz^{n+1}$. See \cite{A3} as well as
\cite{AB3}, \cite{AB2}, and \cite{AH4}, for further references and/or
discussion.

\subsection{Domains in $\sz^2$ and $\hz^2$} \label{geo1}
%Subsection C
Our generalized approach can be used to improve some inequalities
relating the eigenvalues of the Laplace-Beltrami operator on a
bounded domain $\Om$ in $\sz^2$ or $\hz^2$ (with Dirichlet
boundary conditions).  Consider the stereographic projections of
$\sz^n$ to $\rz^n$, for $n\ge 2$, via projection from the south pole of $\sz^n$.
Then the metric is given by (\cite{Cha2}, p.\ 58)
\[ds^2 = p(x)^2 |dx|^2, \quad \text{where} \quad p(x)=\dfrac{2}{1+ |x|^2},\]
where $| \cdot |$ denotes the Euclidean norm.  Hence,
\[g_{ij} = p^2 \delta_{ij}, \quad \mathcal{G} =\left(g_{ij} \right)  =p^2 I, \quad
{\mathcal{G}}^{-1} = \left(g^{ij}\right) = \dfrac{1}{p^2} I, \quad \sqrt g= p^n,\]
where $g=\det \mathcal{G}$ and
%I'm trying to use \det like \sin to represent the determinant.  Perhaps this
%isn't right(?).  It should definitely be checked at some point.
\[\D = \dfrac{1}{p^n} \sum_{i=1}^n
\dfrac{\partial}{\partial x_i} \Big(p^{n-2}
\dfrac{\partial}{\partial x_i}\Big).\]
A Euclidean disk of radius $r$ centered at the origin in $\rz^n$
corresponds to a geodesic disk of radius $\alpha$ in $\sz^n$
centered at the north pole, where $r$ and $\alpha$ are related by
$r = \tan \frac{\alpha}{2}$.

For $n=2$, the Laplace-Beltrami operator takes the form
\[\D_{\sz^2} = \dfrac{1}{p^2} \, \D_{\rz^2}.\]
An eigenvalue of the problem
\[-\D_{\sz^2} u = \la u \quad \text{in} \quad \Om \subset \sz^2,\]
(for $\Om$ a bounded domain) with Dirichlet boundary conditions is
also an eigenvalue of the inhomogeneous membrane problem
\[\D_{\rz^2} u = \la p^2 u\]
also with Dirichlet boundary conditions.  This is then an
inhomogeneous membrane problem with $q(x)=p(x)^2$.  It is obvious
that $q_{max} \le 4.$  Moreover, \[q_{min} = \dfrac{4}{\Big(1+
|x|_{max}^2\Big)^2}.\]
We also have $|x|_{max} = \tan
\dfrac{\Theta}{2}$ by virtue of the correspondence between
geodesic and Euclidean disks, where $\Theta$ is the outer radius
of $\Om$, i.e., the geodesic radius of the circumscribing circle (without
loss of generality, we can assume that this circle is centered on the
north pole).
We have
\[q_{min} = \dfrac{4}{\Big(1+ \tan^2\dfrac{\Theta}{2}\Big)^2}
= 4 \cos^4 \dfrac{\Theta}{2}= \big(1+ \cos \Theta \big)^2.\]
%I played with some \Big's here, changing them to \big's.  These should be
%checked again later.  Similarly in the next displayed inequality.
Therefore,
\[\dfrac{q_{max}}{q_{min}} \leq \dfrac{4}{\big(1+ \cos \Theta \big)^2}.\]
The following is then an extension--\`a la Harrell-Stubbe--of
earlier works by Harrell-Michel \cite{HM2}, Harrell \cite{Ha2},
Cheng \cite{Cheng}, and Ashbaugh \cite{A3}.

\begin{thm}
The eigenvalues of the Laplace-Beltrami operator on a bounded
domain $\Om \subset \sz^2$ with Dirichlet boundary conditions
satisfy the following inequalities
\begin{align}
\sum_{i=1}^m (\la_{m+1} - \la_i)^p \leq  \dfrac{8}{\Big(1+ \cos
\Theta \Big)^2} \, \sum_{i=1}^m (\la_{m+1} - \la_i)^{p-1}  \la_i
\text{ for } p\le 2
%\text{ for } 0 \le p \le 2
\label{HS-sph1}
\end{align}
and
\begin{align}
\sum_{i=1}^m (\la_{m+1} - \la_i)^p \leq  \dfrac{4p}{\Big(1+ \cos
\Theta \Big)^2} \, \sum_{i=1}^m (\la_{m+1} - \la_i)^{p-1}  \la_i
\text{ for } p\ge 2. \label{HS-sph2}
\end{align}
where $0<\Theta<\pi$ designates the outer-radius of $\Om$, i.e., the
radius of the circumscribing geodesic circle.
\end{thm}
\noindent {\em Remark.}  It is to be noted that H.\ C.\ Yang \cite{Y} and
Ashbaugh \cite{A3} produced {\em universal} (i.e., domain
independent) inequalities for $\Om \subset \sz^n$ (see part B of
Section 5 of \cite{A3}).  Following the same arguments one can produce
the following (see Section \ref{immersion} below for more discussion
and the essence of the proof of this theorem).
\begin{thm}
The eigenvalues of the Laplace-Beltrami operator on a bounded
domain $\Om \subset \sz^n$ with Dirichlet boundary conditions
satisfy the following inequalities
\begin{align}
\sum_{i=1}^m (\la_{m+1} - \la_i)^p \leq  \dfrac{1}{n} \,
\sum_{i=1}^m (\la_{m+1} - \la_i)^{p-1} (4 \la_i+n^2) \text{ for }
p\le 2
%\text{ for } 0 \le p \le 2
\label{HS-sph3}
\end{align}
and
\begin{align}
\sum_{i=1}^m (\la_{m+1} - \la_i)^p \leq  \dfrac{p}{2 n}  \,
\sum_{i=1}^m (\la_{m+1} - \la_i)^{p-1} (4 \la_i+n^2)
\text{ for } p\ge 2.
\label{HS-sph4}
\end{align}
\end{thm}

To consider bounds for the eigenvalues of the Laplace-Beltrami operator
on a bounded domain $\Omega \subset \hz^2$, we might consider the
problem using any of several models for $\hz^2$.  We restrict ourselves
to the half-plane model for illustrative purposes.  We refer the reader to
\cite{AH4}, \cite{Ha2}, \cite{HM2} for more discussion.  Here once again
the problem can be thought of as an inhomogeneous membrane problem
(a point of view advocated by Bandle in \cite{Ba})
% This parenthetical remark is likely unnecessary and sort of beside
% the point, as I doubt that Bandle is the person who originated this
% idea, made the original observation, or whatever.
%If we remove this remark, then we should also remove the Bandle
%paper referred to here from our reference list.
since the Laplace-Beltrami operator is given by
\[\D_{\hz^2} = y^2 \D_{\rz^2}.\]

The density function is given by $q({\bf x}) = 1/y^2$ for ${\bf x}
= (x, y) \in \hz^2$.  Our extension then reads.
\begin{thm}
The eigenvalues of the Laplace-Beltrami operator on a bounded
domain $\Om \subset \hz^2$ satisfy the following inequalities
\begin{align}
\sum_{i=1}^m (\la_{m+1} - \la_i)^p \leq  2 \ \dfrac{\sup_{\Om} y^2
}{\inf_{\Om} y^2} \, \sum_{i=1}^m (\la_{m+1} - \la_i)^{p-1},\
\text{ for } p\le 2
%\text{ for } 0 \le p \le 2
\label{HS-hyper1}
\end{align}
and
\begin{align}
\sum_{i=1}^m (\la_{m+1} - \la_i)^p \leq  p \ \dfrac{\sup_{\Om} y^2
}{\inf_{\Om} y^2}  \, \sum_{i=1}^m (\la_{m+1} - \la_i)^{p-1}
\text{ for } p\ge 2.
\label{HS-hyper2}
\end{align}
\end{thm}

\subsection{Eigenvalues of Homogeneous and Minimally Immersed Submanifolds}
%Subsection D
\label{immersion}  Let $M^n$ be an $n$-dimensional compact
manifold (without boundary) of finite volume $V.$  Consider the
problem of estimating the eigenvalues of the Laplace-Beltrami
operator on $M^n$.  The earliest bounds for this problem were
found by Cheng \cite{Cheng} in 1975.  He considered the problem of
estimating these eigenvalues when $M^n$ is immersed in the
Euclidean space $\rz^N$.  Very shortly thereafter, Maeda
\cite{Maeda} considered the analogous problem for domains in the
sphere $\sz^N$ (cf.\ Subsection C above), and for minimally
immersed submanifolds of $\sz^N$.  Also, P.\ C.\ Yang and S.-T.\
Yau \cite{YY} dealt with this problem in the case of a minimally
immersed submanifold of the sphere $\sz^{N}$.  The results of
Maeda and of Yang and Yau (as corrected by Leung
\cite{Leung}) are essentially that
\[\la_{m+1} - \la_{m} \leq n+ \dfrac{2}{n (m+1)} \Big(\sqrt{\Lambda^2+ n^2
\Lambda (m+1)} + \Lambda\Big),\]
where $\Lambda = \sum_{i=1}^m \la_i.$  (We note that $\la_0 = 0$ is the first
eigenvalue for this problem since $M^n$ is compact.)  Beyond that Leung \cite{Leung}, following the approach of Hile and Protter \cite{HP}, produced an HP-type formula in the spirit of Maeda and Yang and Yau.  

In 1995 Harrell and Michel \cite{HM1} (see
also \cite{Michel}, \cite{HM2}) showed, via a general trace inequality, that one
can produce simpler and ``natural'' inequalities which avoid introducing
square root terms such as that found in the bound above.  Finally,
H.\ C.\ Yang \cite{Y} produced, in the same spirit, the strongest version of all
bounds to date.  His 1991 preprint only gradually became known to researchers
in the field.  A revised preprint was circulated in 1995, but neither version was
ever published.  

For further background on the history and context of the methods discussed above one can consult \cite{A3} (see also \cite{AH4}).

Bringing in ideas from Yang \cite{Y} and Harrell-Stubbe \cite{HS}, as developed in this paper, we arrive at the Harrell-Stubbe-type bounds contained in the following theorem.

\begin{thm} \label{thm:minimal}
Let $M^n$ be an $n$-dimensional minimally immersed submanifold of
$\sz^N \subset \rz^{N+1}$, then the eigenvalues of the Laplacian
$-\D_M$, $0=\la_0 \leq \la_1 \leq \la_2 \leq \cdots$, satisfy the
following inequalities
\begin{align}
\sum_{i=0}^m (\la_{m+1} - \la_i)^p \leq  \dfrac{1}{n} \,
\sum_{i=0}^m (\la_{m+1} - \la_i)^{p-1} (4 \la_i+n^2)
\text{ for } p\le 2
%\text{ for } 0 \le p \le 2
\label{HS-min1}
\end{align}
and
\begin{align}
\sum_{i=0}^m (\la_{m+1} - \la_i)^p \leq  \dfrac{p}{2 n} \,
\sum_{i=0}^m (\la_{m+1} - \la_i)^{p-1} (4 \la_i+n^2) \text{ for }
p\ge 2. \label{HS-min2}
\end{align}
\end{thm}
\smallskip \noindent {\em Proof.}
The minimality of the immersion in $\sz^N$ is guaranteed by the
condition that the coordinate functions of the immersion are
eigenfunctions of the Laplace-Beltrami operator on $M^n$ with
eigenvalue $n$. The auxiliary operators are given by the
coordinate functions in this case. Moreover, $\rho_i =n$ and
$\Lambda_i\le n^2 + 4 \la_i$.  Feeding this data into ineqs.\
\eqref{eq:3.142} and \eqref{eq:3.18} (see Corollaries \ref{cor:3.2}
and \ref{cor:3.5}) yields the desired results.
%Decide on the ``right'' Corollary (Corollaries?; inequalities?)
%to refer to here.  Refer to the ineqs.\ in Corollaries 3 and 8.
\hfill $\Box$

Li \cite{Li} dealt with the eigenvalue problem for a compact
homogeneous space. The key to his result and all subsequent
improvements by Harrell and Michel \cite{HM2}, \cite{HM1} (see also
\cite{Michel}, \cite{AH4}) is the following lemma.

\begin{lem} (Li \cite{Li}) \label{lemma:li}
Let $M^n$ be a compact homogeneous manifold of finite volume $V$
and let $\{\phi_{1,\alpha}\}_{\alpha=1}^k$ be a real orthonormal basis
%Did Li insist on the basis eigenfunctions being real?  I would think so,
%though this probably doesn't matter a whole lot.
for the $k-$dimensional eigenspace of the first non-zero
eigenvalue $\la_1$.  Then
\[\sum_{\alpha=1}^{k} \phi_{1, \alpha}^2 = \dfrac{k}{V}
\quad \text{and} \quad \sum_{\alpha=1}^{k} |\Q \phi_{1, \alpha}|^2
\leq \dfrac{\la_1 k}{V}.\]
\end{lem}
%It would appear that the second relation is actually an equality (and is
%a consequence of the first relation, given that we're dealing with eigenfunctions
%of the Laplacian on a compact manifold.  Did Li (and others?) only present
%this as an inequality?

Using this lemma, Li was able to prove that
\[\la_{m+1} - \la_{m} \leq \la_1 + \dfrac{2}{m+1} \Big(\sqrt{\Lambda^2 +
(m+1) \Lambda \la_1} + \Lambda\Big).\] This is of course an
inequality in the spirit of Maeda, Yang-Yau, and Leung (cf.\ also
\cite{A3}).  We have the following improvement (and ``natural
extension'' of the classical inequalities of PPW, HP, and
H.~C.~Yang).
\begin{thm} \label{thm:compact}
Let $M^n$ be a compact homogeneous manifold of finite volume $V$
and let $\{\phi_{1,\alpha}\}_{\alpha=1}^k$ be an orthonormal basis
for the $k-$dimensional eigenspace of the first non-zero
eigenvalue $\la_1$.  Then, its eigenvalues satisfy the following
\begin{align}
\sum_{i=0}^m (\la_{m+1} - \la_i)^p \leq \sum_{i=0}^m (\la_{m+1} -
\la_i)^{p-1} (4 \la_i+\la_1)
\text{ for } p\le 2
%\text{ for } 0 \le p \le 2
\label{HS-cpt1}
\end{align}
and
\begin{align}
\sum_{i=0}^m (\la_{m+1} - \la_i)^p \leq  \dfrac{p}{2} \,
\sum_{i=0}^m (\la_{m+1} - \la_i)^{p-1} (4 \la_i+\la_1) \text{ for
} p\ge 2. \label{HS-cpt2}
\end{align}
\end{thm}
\smallskip \noindent {\em Proof.}
The choices we make are $B_j= \phi_{1,j}$ for $j=1, \cdots, k$
(the eigenfunctions of Lemma \ref{lemma:li}) and $T_j=[-\D,
B_j].$  Then, $\rho_i =\dfrac{\la_1 k}{V}$ and $\Lambda_i \leq
\dfrac{\la_1 k}{V} \ \big(4  \la_i+ \la_1\big).$  Hence the desired
results follow via Corollaries \ref{cor:3.2} and \ref{cor:3.5} (ineqs.\
\eqref{eq:3.142} and \eqref{eq:3.18}).
%Again, we need to decide on the ``right'' Corollary (Corollaries?;
%inequalities?) to refer to here.  Refer to the ineqs.\ in Corollaries
%3 and 8.
\hfill $\Box$
%In this Li business, one should really be able to bring back an $n$
%in the denominator.  I don't know why no one has really worked on
%that.  My impression is that Li was unenterprising to begin with, and
%then everyone else just followed him.

\subsection{Second Order Elliptic Operators} \label{elliptic}
%Subsection E
In \cite{Ho2}, Hook considered a general, second order, elliptic
partial differential equation with constant coefficients of the form
\begin{align}
\mathcal{A} \, u \equiv -\sum_{i,j=1}^n \dfrac{\partial}{\partial x_i}
\big(a_{ij} \dfrac{\partial u}{\partial x_j}) + \sum_{i=1}^n b_i
\dfrac{\partial u}{\partial x_i} = \la \, u \label{eq:elliptic.1}
\end{align}
on a bounded domain $\Om \subset \rz^n$ with Dirichlet boundary
conditions.  With the assumption that $A = [a_{ij}]$ is a symmetric
positive definite real matrix, he was able to produce HP-type
bounds for the eigenvalues of this problem.  In \cite{AH4} we
succeeded in producing H.~C.~Yang-type bounds for this problem
thus strengthening Hook's results.  The essential ingredient is to
rewrite the problem in the form
\begin{align}
- e^{w \cdot x} \Div(A \,  e^{-w \cdot x} \, \Grad \ u) = \la u,
\label{eq:elliptic.4}
\end{align}
where $w \in \rz^n$ is a constant vector given by $w=A^{-1} b$,
and $b=[b_i]$ appearing in equation (\ref{eq:elliptic.1}).  The
matrix $A$ is diagonalized according to $A=U^{-1} K U$, with
$U$ a real orthogonal matrix.  The standard basis $e_1, e_2,
\cdots, e_n$ is then transformed according to $v_j=U^{-1} e_j$ to
produce a new orthonormal basis $v_1, v_2, \dots, v_n.$  The
operators $T_j$ are given by
\begin{align} T_j \, u
= (v_j , \sqrt A \ \Grad \ u) - \dfrac{1}{2} (v_j , \sqrt A \, w) \, u,
\label{eq:elliptic.5}
\end{align}
where $( \cdot , \cdot )$ denotes the usual dot product in
$\rz^n,$ $\sqrt A$ denotes the positive definite square root of $A,$
and $j=1, 2, \dots, n.$  The operators $T_j$ are skew-symmetric
with respect to the inner product
\begin{align}
\langle u , v \rangle = \int_{\Om}  u \bar{v} e^{-w \cdot x} dx.
\label{eq:elliptic.6}
\end{align}
Also, our original operator $\mathcal{A}$ satisfies
\begin{align}
\langle \mathcal{A} \, u, u \rangle = \sum_{j=1}^n \langle T_j u,
T_j u \rangle + \dfrac{1}{4} ( A w , w)  \,  \langle u , u \rangle
\label{eq:elliptic.7}
\end{align}
in the given inner product.  The auxiliary operators $B_k$ are
chosen to be of the form $B u = \phi(x) u$, a multiplication by a
real-valued function $\phi$ of the coordinates.  The commutation
conditions $[T_j, B_k] = \delta_{j k}$ are equivalent to
\begin{align}
(\sqrt A \ v_j, \Grad \ \phi_k) = \delta_{jk}.
\label{eq:elliptic.11}
\end{align}
The vectors $\{v_1, v_2, \dots, v_n\}$ form a basis for $\rz^n$
and the same is the case for $\{\sqrt A \ v_1, \sqrt A \ v_2,
\dots,\sqrt A \ v_n\}$ since $\sqrt A$ is invertible.  We form the
matrix $C$ with columns given by the elements of $\sqrt A \ v_1,$
$\sqrt A \ v_2,$ $\dots,$ $\sqrt A \ v_n.$  C is then invertible.
We let $F=[f_{jk}]$ be its inverse.  Condition
(\ref{eq:elliptic.11}) is equivalent to \[\Phi \, C = I,\] where
$\Phi$ is the matrix with rows given by \[\Grad \ \phi_1, \Grad
\ \phi_2, \dots, \Grad \ \phi_n,\] and $I$ is the identity matrix.
Hence, $F=C^{-1}=\Phi$ and
\[\dfrac{\partial \phi_j}{\partial x_k} = f_{jk}.\]
The functions \[\phi_j = \sum_{j=1}^n f_{jk} \ x_k\] satisfy the
conditions we seek.
\begin{thm}
With $M= \sqrt A$, the eigenvalues of problem
(\ref{eq:elliptic.1}) satisfy the inequalities
\begin{align}
\sum_{i=0}^m (\la_{m+1} - \la_i)^p \leq  \dfrac{4}{n} \,
\sum_{i=0}^m (\la_{m+1} - \la_i)^{p-1} \big(\la_i - \dfrac{1}{4} \
\|M^{-1} b \|^2\big) \text{ for } p\le 2 \label{HS-ellip1}
\end{align}
and
\begin{align}
\sum_{i=0}^m (\la_{m+1} - \la_i)^p \leq  \dfrac{2p}{n} \,
\sum_{i=0}^m (\la_{m+1} - \la_i)^{p-1} \big(\la_i - \dfrac{1}{4} \
\|M^{-1} b \|^2\big) \text{ for } p\ge 2. \label{HS-ellip2}
\end{align}
\end{thm}
%Shouldn't there be a proof, or a reference to other results that make
%the proof obvious here?  Presumably one should at least make a
%connection to Schroedinger operators H=A+V, with A given as the
%negative of the sum of the Ts squared, or some such here.  On the
%other hand, by eq. (83) (bottom of the previous page) A is not exactly
%of this form, though it only differs from that form by a constant (which
%can be absorbed into the eigenvalue parameter).  Or one could
%follow Hook in defining the operator A with this constant already
%subtracted out.
%An alternative proof, which gets one this result directly, is outlined
%in Remark 2 below.
\noindent{\em Remarks.}  1.  The use of the orthonormal vectors
$v_j$ in our definition of the operators $T_j$ in
\eqref{eq:elliptic.5} is not necessary either for the
skew-symmetry of the $T_j$'s in the inner product $\langle \cdot,
\cdot \rangle$ nor for the identity \eqref{eq:elliptic.7}; for
both of these it is enough that the vectors $\{v_j\}_{j=1}^n$ form
an orthonormal basis.

\noindent  2.  A more direct approach to this result, since here
we deal only with the case of constant coefficients, is simply to
transform away the symmetric matrix $A$ via a linear change of
variables (in fact, the change from $x_k$ to $\phi_j$ given
above), arriving at a transformed problem on a new bounded domain
$\tilde{\Omega}$ in the variables $\phi_j$ where the second-order
part of the operator is just the Laplacian (in the variables
$\phi_j$, with the $\phi_j$'s viewed as Euclidean variables) and
with a first-order part involving a new $b$-vector,
$\tilde{b}=M^{-1}\,b$.  The first-order term can then be entirely
eliminated via the change of dependent variable,
$v=e^{-(\tilde{b},\phi)/2}\,u$ (here $\phi$ denotes the vector
having the $\phi_j$'s as components), producing an eigenvalue
problem $-\Delta \, v = \mu \, v$ for the Laplacian on a bounded
domain, still with homogeneous Dirichlet boundary conditions, and
with the only modification being that the eigenvalue parameter
$\lambda$ becomes $\mu = \lambda-\|\tilde{b}\|^2/4$.  Finally,
since we already know the inequalities \eqref{eq:intro.6} and
\eqref{eq:intro.5} for the Laplacian, we obtain the results of the
theorem simply by replacing all $\lambda$'s in those inequalities
by $\mu$'s, where $\mu_i =\lambda_i - \|M^{-1}\,b\|^2/4$.

\subsection{Sturm-Liouville Problem}
\label{sturm}
%Subsection F
Let $I=(a, b) \subset \rz$.  Hook \cite{Ho2} considered the
following Sturm-Liouville problem on $I$
\begin{align}
\mathcal{A}\,u &= - (p \, u')' + q u = \la u, \notag \\
u(a) &= u(b) = 0, \label{eq:sturm.1}
\end{align}
where $p(x) >0$  and $q(x)$ are real-valued functions on $I.$  The
differential operator $\mathcal{A}$ is symmetric with respect to
the inner product $\langle u, v \rangle = \int_a^b u \bar{v}\, dx.$
One is able to prove that (see \cite{AH4}, \cite{Ho2})
\begin{align}
\langle \mathcal{A} \, u, u \rangle = \langle T \, u, T \, u \rangle
+ \langle \mathcal{Q} \, u, u \rangle.
\end{align}
where $T \, u = \dfrac{1}{2} \big( \sqrt p \ u' + (\sqrt p \
u)'\big),$ and $\mathcal{Q}\, u = Q(x) \, u$ for \[Q(x) = q(x) -
\dfrac{1}{16} \dfrac{p'(x)^2}{p(x)} +\dfrac{1}{4} p''(x).\] $T$ is
skew-symmetric.  The symmetric operator $B$ is chosen to be of the
form $B \, u = \phi(x) \, u$ with $\phi$ real-valued.  The commutation
condition $[T, B]= 1$ yields the form \[\phi(x) = \int \dfrac{dx}{\sqrt {p(x)}}.\]
The following theorem is then immediate.

\begin{thm} Suppose $Q(x) \ge M$ for $x\in (a,b)$. Then,
the eigenvalues of problem (\ref{eq:sturm.1}) satisfy the
inequalities
\begin{align}
\sum_{i=0}^m (\la_{m+1} - \la_i)^p \leq 4 \, \sum_{i=0}^m
(\la_{m+1} - \la_i)^{p-1} \big(\la_i - M\big) \text{ for } p\le 2
\label{HS-sturm1}
\end{align}
and
\begin{align}
\sum_{i=0}^m (\la_{m+1} - \la_i)^p \leq  2 p \, \sum_{i=0}^m
(\la_{m+1} - \la_i)^{p-1}  \big(\la_i - M\big) \text{ for } p\ge
2. \label{HS-sturm2}
\end{align}
\end{thm}
\noindent {\em Remark.} Hook \cite{Ho2} considered two more
problems for which he produced HP-type bounds. The first is a
diagonal $n$-dimensional version of the Sturm-Liouville problem.
The second is a Schr\"odinger operator with magnetic potential.
For both problems H.\ C.\ Yang-type inequalities were produced in
\cite{AH4} and as such extensions \`a la Harrell-Stubbe are valid
as well.

\begin{recrev}
Received 28 December 2007\\
Revised xxx 200x
\end{recrev}

\end{document}